\documentclass{ejm}
\usepackage{graphicx}
\usepackage{hyperref}
\usepackage{amsmath}
\newtheorem{thm}{Theorem}[section]

\newtheorem{lem}[thm]{Lemma}
\newtheorem{prop}[thm]{Proposition}
\newtheorem{rem}[thm]{Remark}

\let\realverbatim=\verbatim
\let\realendverbatim=\endverbatim
\renewcommand\verbatim{\par\addvspace{6pt plus 2pt minus 1pt}\realverbatim}
\renewcommand\endverbatim{\realendverbatim\addvspace{6pt plus 2pt minus 1pt}}

\ifprodtf \else
\checkfont{eurm10}
\iffontfound
\IfFileExists{upmath.sty}
{\typeout{^^JFound AMS Euler Roman fonts on the system,
		using the 'upmath' package.^^J}%
	\usepackage{upmath}}
{\typeout{^^JFound AMS Euler Roman fonts on the system, but you
		don't seem to have the}%
	\typeout{'upmath' package installed. EJM.cls can take advantage
		of these fonts,^^Jif you use 'upmath' package.^^J}%
	\providecommand\mathup[1]{##1}%
	\providecommand\mathbup[1]{##1}%
}
\else
\providecommand\mathup[1]{#1}%
\providecommand\mathbup[1]{#1}%
\fi
\fi

\ifprodtf \else
\checkfont{msam10}
\iffontfound
\IfFileExists{amssymb.sty}
{\typeout{^^JFound AMS Symbol fonts on the system, using the
		'amssymb' package.^^J}%
	\usepackage{amssymb}%
	  \let\leq=\leqslant
	  \let\geq=\geqslant
}{}
\fi
\fi

\ifprodtf \else
\IfFileExists{amsbsy.sty}
{\typeout{^^JFound the 'amsbsy' package on the system, using it.^^J}%
	\usepackage{amsbsy}}
{\providecommand\boldsymbol[1]{\mbox{\boldmath $##1$}}}
\fi

\newsavebox{\astrutbox}
\sbox{\astrutbox}{\rule[-5pt]{0pt}{20pt}}

\usepackage{xcolor}
\usepackage{soul}

\newdefinition{definition}[theorem]{Definition}
\newtheorem{lemma}{Lemma}

\title[European Journal of Applied Mathematics]{Nonlinear Biphasic Mixture Model: Existence and Uniqueness Results}

\author[M. Alam et al.]{%
	M.\ns A\ls L\ls A\ls M$\,^1$,\ns
	A.\ns  M\ls U\ls N\ls T\ls E\ls A\ls N$\,^2$\ns
	\and
	G.\ns P.\ns R\ls A\ls J\ls A\ns S\ls E\ls K\ls H\ls A\ls R$\,^3$
}

\affiliation{%
	$^1\,$Department of Mathematics, Ecole Centrale  School of Engineering, Mahindra University, Hyderabad, 500043, Telangana, India\\
	$^2\,$Department of Mathematics and Computer Science, University of Karlstad,  Universitetsgatan 2, 651 88 Karlstad, Sweden\\
	$^3\,$Department of Mathematics, Indian Institute of Technology Kharagpur, WB-721302, India\\
	email\textup{\nocorr: \texttt{rajas@iitkgp.ac.in}}}
\date{10 February 2022}
\pubyear{2000}
\volume{000}
\pagerange{\pageref{firstpage}--\pageref{lastpage}}

\begin{document}
	
	\label{firstpage}
	\maketitle
	
	\begin{abstract}
		This paper is concerned with the development and analysis of a mathematical model that is motivated by interstitial hydrodynamics and tissue deformation mechanics
		(poro-elasto-hydrodynamics) within an in-vitro solid tumor. The
		classical mixture theory is adopted for mass and momentum balance
		equations for a two-phase system. A main
		contribution of this study, we treat the physiological transport parameter (i.e., hydraulic resistivity) as anisotropic and heterogeneous, thus the governing system is strongly coupled and non-linear. 
		We derived a weak formulation and then formulated the equivalent fixed-point problem. This enabled us to use the Galerkin method, and the classical results on monotone
		operators combined with the well-known Schauder and Banach {fixed point} theorems to prove the existence and uniqueness results. 
	\end{abstract}
	
	
	\begin{subjclass}[2020]
		76Txx, 76Zxx (Primary); 
		35Q74, 35D30 (Secondary)
	\end{subjclass}
	
	\section{Introduction}
	
	The study of fluid flow through porous media has gained the
	attention of many researchers over the years. Examples of natural porous materials are
	living tissues, rocks, soils, etc. On the other hand, there are manmade porous materials, e.g., concrete, foam, ceramic, etc.
	Because of their wide presence and hard-to-estimate effective properties, flow through porous material is studied by engineers and scientists.
	This leads to a coupled phenomenon where fluid flow and solid deformation in porous materials interplay. This is a
	classical problem in geomechanics and biomechanics. Recently, one of the most studied topics in the field of fluid mechanics is flow through biological tissues such as
	tumors, Gylcocalyx layers, and articular cartilage.
	This paper deals with the mathematical modeling and analysis of
	the coupled phenomena of fluid flow and solid deformation (in
	short \textit{poro-elasto-hydrodynamics}) within an in-vitro tumor model.
	Typically, a tumor is assumed to be a deformable porous medium
	that consists of multiple phases, e.g., one fluid phase and many solid ones. A tumor may exist in isolation (in-vitro)  or may be present in normal tissue (in-vivo) \cite{dey2016hydrodynamics}. Although the internal geometrical structure of tumors is complicated,
	developing mathematical models for approximate situations is very
	useful. Theoretical predictions generated from such approximate
	models may help to reduce the number of animal experiments
	that need to be carried out and also suggest new experimental
	programmes that identify optimal tumor therapy schedules
	\cite{byrne1999using}.
	
	Mathematical models of tumors in general   can be divided into three categories:
	discrete, continuum, and hybrid. Here, we focus on continuum models
	that treat cells as averaged populations and are based on the
	continuum mechanics approach to porous media combined with mixture theory
	\cite{shelton2011mechanistic}.
	The early mathematical models on
	the avascular tumor growth \cite{araujo2004history} assumed that
	tumors are made of single
	type of cells having a constant density.  However, various experimental, and theoretical	evidence have shown that such a description is not sufficient to study the tumor dynamics
	\cite{araujo2004history,preziosi2009multiphase}. Hence,  multiphase models came into play.
	In this case, one can consider density variations within mixture components to evaluate the evolution of partial stresses. 
	Biot's theory of poroelasticity and the theory of mixture are
	commonly adapted models to explore
	\textit{poro-elasto-hydrodynamics} \cite{dey2016hydrodynamics}.
	Alike continuum level approaches involve the development of a set of
	equations to represent the mechanical behaviour of a soft tissue
	such as a tumor (which is assumed as a deformable porous material)
	at the macro scale, using a porous media approach.
	
	The first multiphase model for tumors is
	proposed by Please et al. \cite{please1998new}. The authors
	proposed a diffusion equation for cell concentration and
	generalized Darcy equation for cells and water motion inside the
	tumor.
	Further, such multiphase models have been studied analytically and numerically by several authors, one can refer \cite{ambrosi2002closure,byrne2003modelling,preziosi2009multiphase}.
	Sumets et al. \cite{sumets2015boundary} described a new boundary-integral representation for biphasic
	mixture theory, where they solved elastohydrodynamic–mobility problems using boundary element methods.
	Dey and Raja Sekhar \cite{dey2016hydrodynamics} used a biphasic mixture model to poro-elasto-hydrodynamics and nutrient transport inside an in-vitro solid tumor. The authors assumed the presence of an unknown fluid source/sink in the model and
	biphasic mixture equations have been solved explicitly in the case of one-dimensional spherical geometry.	
	Slvia and Wheeler \cite{barbeiro2010priori} presented a coupled geomechanics and unsteady reservoir flow model
	using the theory of poroelasticity. They established the existence and uniqueness of a weak solution and computed \textit{a priori} error
	estimates for the numerical solution with stress-dependent permeability. {In \mbox{\cite{jager2011homogenization}}, a nonlinear model for a poroelastic medium (described by quasi-static Biot-equations) coupled to transport equations of substances was considered. They have modeled time and space-dependent processes in deformable cellular tissues by the method of homogenization, starting from a reactive flow system coupling mechanics at the pore scale. The model was analyzed and the global-in-time existence and uniqueness of the solution were shown.}
	Cao et al. \cite{cao2014steady} have considered
	a nonlinear steady flow model in a deformable biological domain based on the theory of poroelasticity (nonlinearity is
	due to the assumption of dilation-dependent interstitial permeability of the solid matrix). They established the existence
	and uniqueness of a weak solution. 
	Looking through mentioned literature, we observe that there is a gap in dealing with the existence and uniqueness of corresponding biphasic mixture equations that describe the coupled phenomena of fluid flow and solid deformation within  biomaterials such as tumors.
	Attempting to fill such a gap, Alam et al.
	\cite{alam2018mathematical,alam2019mathematical} developed a
	well-posedness theory and certain regularity results in $2d$ and $3d$ for poro-elasto-hydrodynamics
	model inside an in-vitro solid tumor. Further, in the case of an in-vivo solid tumor, Alam et al.  \cite{alam2021existence} developed existence and uniqueness results in a weak sense for poro-elasto-hydrodynamics while assuming the	hydraulic resistivity heterogeneous and deformation dependent.

	We note that, in general,  poro-elasto-hydrodynamics models
	within a tumor may not lead directly to linear biphasic
	mixture equations. In practice, due to the
	non-uniform blood vessel distribution, the supply of fluids and macromolecules within a tumor is heterogeneous. As
	a consequence, physiological transport parameters (e.g., hydraulic
	resistivity or permeability) depend on space and deformation. {For
		instance, in the case of soft permeable tissue and gel,} Barry and Aldis
	\cite{barry1990comparison}, Holmes and Mow
	\cite{holmes1990nonlinear} {considered permeability depending
		exponentially on the strain.} Also, some of the biological tissues
	and cells display anisotropic permeability
	\cite{ateshian2010anisotropic}. In particular, articular cartilage
	typically exhibits anisotropic behavior
	\cite{federico2008anisotropy}. {Further, in the case of multicellular
		tumor tissue}, Giverso and
	Preziosi \cite{giverso2019influence} follow Holmes and Mow
	\cite{holmes1990nonlinear} and {considered that the permeability depends on the volumetric
		deformation (or strain).} 
	Dey and Sekhar
	{\cite{dey2016hydrodynamics}} {assumed that the hydraulic conductivity of soft tumor depends on the radial distance.}
	These nonlinear effects inserted in the physiological parameter yield nonlinear biphasic mixture equations.
	As far as we know, for nonlinear biphasic mixture models, there is a lack of literature regarding the existence, uniqueness, and regularity of the solution. In this paper, we present a nonlinear biphasic mixture model that represents poro-elasto-hydrodynamics which do not account for any new growth of tumor cells.
	The physiological transport parameter (hydraulic resistivity) is
	assumed to be deformation dependent, which yields the nonlinearity in the model. We develop a local weak solvability theory.
	
	\subsection{Biphasic Mixture Theory}\label{Bi}
	In this subsection, we introduce the generic governing equations. We use biphasic mixture theory to represent
	the fluid and solid phases of the tumor. Following \cite{barry1991fluid,byrne2003modelling,dey2016hydrodynamics}, we
	apply the conservation of mass and momentum to the fluid and solid phases, viewing the fluid as viscous Newtonian and the
	solid as deformable, and accounting for momentum exchange between the two phases. Let $\mathbf{V}_f$ and $\mathbf{V}_s$
	denote the velocities of the fluid and solid phases, respectively. The apparent densities of the fluid and solid phases
	are denoted   by $\tilde{\rho}_f$ and $\tilde{\rho}_s,$    respectively, and their corresponding volume fractions by
	$\phi_f$ and $\phi_s.$ The true densities of the fluid and solid phases are then $\rho_f=\phi_f\tilde{\rho}_f$ and
	$\rho_s=\phi_s\tilde{\rho}_s.$     Accordingly, in $\Omega$ the mass and linear momentum balance equations
	for the fluid phase are given by
	\begin{equation}\label{Eq1}
		\frac{\partial(\tilde{\rho}_f\phi_f)}{\partial t}
		+\nabla\cdot[(\tilde{\rho}_f\phi_f)\mathbf{V}_f]= \tilde{\rho}_f
		S_f,
	\end{equation}
	\begin{equation}\label{Eq11}
		\rho_f\left(\frac{\partial\mathbf{V}_f}{\partial
			t}+(\mathbf{V}_f\cdot\nabla)\mathbf{V}_f \right)= \nabla\cdot
		\mathbf{T}_f+\boldsymbol{\Pi}_f+\mathbf{b}_f,
	\end{equation}
	where $\textbf{T}_{f}$  denotes the stress tensor for the fluid
	phase
	\begin{eqnarray}
		\mathbf{T}_{f}=-[\phi_{f}P-\lambda_{f}\nabla\cdot\mathbf{V}_{f}]\mathbf{I}
		+\mu_{f}[\nabla \mathbf{V}_{f}+(\nabla \mathbf{V}_{f})^{T}].
	\end{eqnarray}
	The corresponding  mass and linear momentum equations for the solid phase are
	\begin{equation}\label{Eq2}
		\frac{\partial(\tilde{\rho}_s\phi_s)}{\partial t}
		+\nabla\cdot[(\tilde{\rho}_s\phi_s)\mathbf{V}_s]= \tilde{\rho}_f
		S_s,
	\end{equation}
	\begin{equation}\label{Eq11a}
		\rho_s\left(\frac{\partial\mathbf{V}_s}{\partial
			t}+(\mathbf{V}_s\cdot\nabla)\mathbf{V}_s \right)=
		\nabla\cdot\mathbf{T}_s+\boldsymbol{\Pi}_s+\mathbf{b}_s,
	\end{equation}
	where $\mathbf{T}_{s}$  denotes the stress tensor for  the solid
	phase
	\begin{eqnarray}\label{Eq12}
		\mathbf{T}_{s}=-[\phi_{s}P-\chi_s
		(\nabla\cdot\textbf{U}_{s})]\mathbf{I}+\mu_{s}[\nabla
		\mathbf{U}_{s}+(\nabla \mathbf{U}_{s})^{T}].
	\end{eqnarray}
	In equations (\ref{Eq1}) and (\ref{Eq2}),  $S_f$ and $S_s$ are fluid
	and solid source terms respectively. $\mathbf{U}_s$ denotes the
	displacement of the solid phase. Hence
	$\mathbf{V}_{s}=\frac{\partial\mathbf{U}_{s}}{\partial t}.$  The average interstitial fluid pressure (IFP) is \(P\) and
	$\mathbf{b}_j$ $j=\{1,2\}$ denotes the body force. Further, the volume
	fractions $\phi_f$ and $\phi_s$ are assumed to satisfy the
	following saturation assumption
	\begin{equation}\label{Eq3}
		\phi_f+\phi_s= 1.
	\end{equation}
	We suppose that the two phases interact together via the drag forces	$\boldsymbol{\Pi}_{s}$ and $\boldsymbol{\Pi}_{f},$ which by Newton's third law, are equal and opposite. Following
	\cite{ambrosi2002closure,byrne2003modelling}, we define
	\begin{eqnarray}\label{Eq12a}
		-\boldsymbol{\Pi}_{s}=\boldsymbol{\Pi}_{f}=
		{K}(\textbf{V}_{s}-\textbf{V}_{f})- (\nabla\phi_{s})P,
	\end{eqnarray}
	where $K=\frac{\mu_f}{k}$ is the hydraulic resistivity or drag coefficient, where $k$ is the permeability of the porous
	matrix (the precise nature of $K$ will be defined in the next section).  Furthermore, $\mu_f$ ($\mu_s$) is the dynamic viscosity
	of the fluid phase
	(solid phase), while $\lambda_f$ and $\chi_s$ denote the Lame
	coefficient and shear modulus of the fluid and solid phases. The
	elastic modulli $\chi_s$ and $\mu_s$ are related to the Young's
	modulus $(\mathcal{Y})$ and Poisson's ratio $(\nu_p)$ via the relationship
	$$\chi_s=\frac{\nu_p\mathcal{Y}}{(1+\nu_p)(1-2\nu_p)}~\mbox{and}~\mu_s
	=\frac{\mathcal{Y}}{2(1+\nu_p)}.$$
	\subsection{Main modeling assumptions}\label{assump}
	Having presented a more generic set of mixture theory equations in
	the absence of any assumptions or boundary conditions, here we
	list some biologically suitable assumptions restricted to a
	specific tumor model. Our choice of model is motivated by the study of an isolated (in-vitro) tumor that behaves as a heterogeneous deformable porous medium. Usually, tumor tissues are considered as incompressible fluids with no voids present. It is assumed that each phase has an equal constant density.
	Suppose $\Omega\subset\mathbb{R}^d,$ $d=\{2,3\}$ is a bounded Lipschitz domain that is filled by the tumor. Let $\partial\Omega$ be its boundary.
	\begin{figure}
		\centering
		\includegraphics[width=6cm]{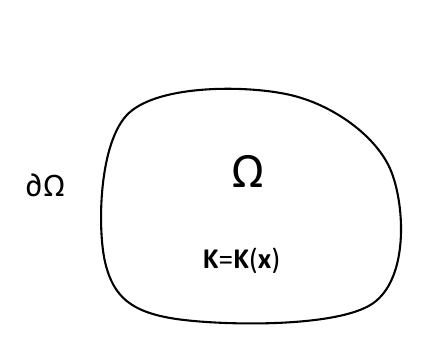}
		\caption{Geometry of the problem} \label{chap2F1}
	\end{figure}
	One may note that the solid tumor is essentially a multicellular
	spheroid. When nutrients perfuse the interstitial space, a large
	number of cells receive adequate food for survival and
	proliferation. As a consequence, the tumor grows in size. For a
	growing tumor, the permeability and the effective mechanics parameters ($\chi_s$
	and $\mu_s$) may depend on the volume fraction of the cell
	population {\cite{byrne2003modelling}}. Moreover, the volume fractions
	depend on space. Hence, it is extremely difficult to analyze mathematically the
	growth model and fluid transport model simultaneously. 
	To simplify this, we assume that the tumor tissue is not
	growing and all elastic parameters ($\chi_s$ and $\mu_s$ etc.) are
	independent of volume fractions.
	Further, we make the following modeling assumptions ({\bf A1})-({\bf A5}):
	\begin{description}
		\item[(A1)]
		Nutrient perfusion and transport occur on much shorter timescales
		than the timescale for tumor cell growth. Accordingly, we view the tumor as a static, perfused biological domain. On the short timescale associated with nutrient transport to (and within) the tumor, cell death, and proliferation are assumed to be negligible. Therefore, we	fix $S_s = 0$ in the tumor and normal tissue regions. Further, on
		the timescale of interest, the solid volume fraction $\phi_{s}$
		remains constant, and for simplicity, we assume it to be
		independent of spatial position and time as (see
		\cite{dey2016hydrodynamics}).
		\item[(A2)] A fluid source is attached to the mixture, hence $S_f\neq0,$ \cite{dey2016hydrodynamics,netti1997macro}.
		\item[(A3)] Motion of the cells and interstitial fluid flow is so slow that the inertial terms can be neglected in both phases, see e.g. \cite{wang1995effect,byrne2003modelling}.
		\item[(A4)] {\bf Structure of the hydraulic resistivity:}
		Various
		experimental and theoretical investigations indicate clearly that
		for the deformable porous medium (or soft biological tissue such
		as articular cartilage, arterial tissue, and tumor), the permeability
		also called in this context hydraulic resistivity depends on stress, dilatation, volume
		fractions (porosity), etc. There are several analytical expressions for
		permeability that are used in literature such as $k(x)=\exp[mU'(x)],$
		where $m$ is a constant and $U'=\frac{dU}{dx}$ with $U$ as
		displacement, while $k(x)=k_0/[1 -
		mU'(x)]$ or $k(x)=k_0[1 + mU'(x)]$ for small $mU'(x)$, or $k(x)=k_0(\phi/\phi_0)^n$ \cite{barry1990comparison,holmes1990nonlinear}.
		Here $k_0$ is the permeability at reference porosity $\phi_0$ and
		$n$ is a variable that may be determined by fitting experimental
		data. In this framework, we propose two different cases. We assume that the
		hydraulic resistivity $K$ depends explicitly on (a)   the solid
		phase displacement $\mathbf{U}_s$ and (b) the strain/dilatation
		i.e. on $\nabla\cdot\mathbf{U}_s.$ 
		In both of these cases, we admit
		anisotropic effects, i.e. $\mathbf{K}$ is a square
		matrix of order $d =\{2,3\}.$ {Note that in case (a) we do
			not use any specific, explicit expression of
			$\mathbf{K}(\mathbf{U}_s)$ in our analysis. One may think though
			of the following choice}
		$${\mathbf{K}(\mathbf{U}_s)=(\alpha_1|\mathbf{U}_s|+\alpha_3)\mathbf{I}
			+(\alpha_1-\alpha_2)
			\frac{\mathbf{U}_s\mathbf{U}_s^T}{|\mathbf{U}_s|},}$$ {where
			$\alpha_i$ ($i=1,~2,~3$) are real constants such that
			$\alpha_1\geq0,\alpha_2\geq0$ and $\alpha_3>0.$ One can easily
			show that $\mathbf{K}$ is Lipschitz and uniformly positive
			\mbox{\cite{sun2002combined}}.}
		{For case (b), one can think of a form  $\mathbf{K}(\nabla\cdot\mathbf{U}_s) =
		 (\gamma_1+\gamma_2|\nabla\cdot\mathbf{U}_s|)\mathbf{I},$ where
		 $\gamma_i\geq 0$ are real constants, and $\mathbf{I}$ is the
		 identity matrix. One can observe
			 that $\mathbf{K}$ is Lipschitz and uniformly positive. 
  We will state further assumptions on \({\bf K}\) in the next Section \ref{s1},
  \cite{barry1990comparison,holmes1990nonlinear}.}
	\end{description}
	\noindent The assumptions ({\bf A1})-({\bf A4}) are taken into account
	in equations (\ref{Eq1})-(\ref{Eq12a}). This helps us to define a coupled nonlinear
	system  of steady-state mass and momentum balance equations in unknowns \((\mathbf{V}_f,\mathbf{U}_s, P)\) as follows:
	\begin{equation}\label{Eq15}
		-\nabla\cdot(2\mu_fD(\mathbf{V}_f)+\lambda_f(\nabla\cdot\mathbf{V}_f)\mathbf{I}
		-\phi_fP\mathbf{I})
		+\mathbf{K}(\boldsymbol{\varsigma})\mathbf{V}_f=\mathbf{b}_f~~\mbox{in}~~\Omega,
	\end{equation}
	\begin{equation}\label{Eq16}
		-\nabla\cdot(2\mu_sD(\mathbf{U}_s)+\chi_s(\nabla\cdot\mathbf{U}_s)\mathbf{I}-\phi_sP\mathbf{I})
		-\mathbf{K}(\boldsymbol{\varsigma})\mathbf{V}_f=\mathbf{b}_s~~\mbox{in}~~\Omega,
	\end{equation}
	\begin{equation}\label{Eq17}
		\nabla\cdot(\phi_f\mathbf{V}_f)=S_f~~\mbox{in}~~\Omega,
	\end{equation}
	where $\boldsymbol{\varsigma}$ is either $\nabla\cdot\mathbf{U}_s$ or
	$\mathbf{U}_s.$ Further, $D(\cdot)$ denotes the
	deviatoric matrix which is defined as
	$D(\mathbf{u})=\frac{1}{2}(\nabla\mathbf{u}+(\nabla\mathbf{u})^t),$ where $(\nabla\mathbf{u})^t$ denotes the transpose of the matrix
	$\nabla\mathbf{u}.$ We have made the approximation
	$\textbf{V}_{s}=\frac{\partial\mathbf{U}_{s}}{\partial t}\approx0$
	which is consistent with the infinitesimal strain theory that is used in e.g. \cite{damiano1996axisymmetric,dey2016hydrodynamics}. The mass
	balance equations for the fluid phase include a source term
	$S_f$, which models fluid exchange with vasculature and lymph
	vessels. For a closed mixture (for instance the case of avascular tumors), we
	consider $S_f=0$ so that $\nabla\cdot(\phi_f\mathbf{V}_f)=0.$
	This is the counterpart of the incompressibility
	constraint.
	Note that in equation (\ref{Eq17}), even though the density of each phase is constant, the vector $\mathbf{V}_f$ is not solenoidal. On the other hand, when the external sources/sinks are attached to the mixture, we have $S_f\neq0$ \cite{netti1997macro,preziosi2009multiphase}.
	Typically, the fluid source $S_f$ is assumed to be driven by the average transmural pressure and (trusting \cite{netti1997macro,dey2016hydrodynamics}) takes the form
	\begin{equation}\label{Eq9a}
		S_f=-{L_p}\left(\frac{A}{V}\right)\{1+L_rA_r\}(P-P_F),
	\end{equation}
	where  $L_p$ is the average hydraulic conductivity coefficient of
	capillary. In (\ref{Eq9a}) $A/{V}$ denotes the capillary surface area per unit
	tissue volume in the tumor tissue and $L_rA_r$ denotes the ratio of
	the strength of distributed solute source through the vasculature
	and solute sink through the lymph vessels and $P_F$ is the
	weighted vascular pressure.
	
	\subsection{{Boundary conditions}}
	The model that we have considered here is supposed to mimic an in-vitro tumor.
	Accordingly, we
	prescribe
	\begin{equation}\label{Eq18}
		\mathbf{T}_f\cdot\mathbf{n}=\mathbf{T}_\infty~~\mbox{and}~\mathbf{U}_s=0~\mbox{on}~~\partial\Omega,
	\end{equation}
	{where $\mathbf{n}$ is the outward normal unit vector to the boundary $\partial\Omega.$} 
	\section{{Non-dimensional equations}}
	Using the transformations
	$\hat{\textbf{x}}=\frac{\textbf{x}}{R},$ $\nabla'=R\nabla,$  $\hat{P}=\frac{P}{P_F},$ $\hat{\mathbf{V}}_f
	=\frac{\mathbf{V}_f}{\frac{RP_F}{\mu_f}},$ $\hat{\mathbf{U}}_s=\frac{\mathbf{U}_s}{\frac{R^3P_F}{\mu_f\nu}},$ 
	$\mathbf{K}=\hat{\mathbf{K}}K_d,$ where $K_{d}$ is the hydraulic
	resistivity (drag coefficient) of the tumor tissue in the absence
	of deformation, $R$ is the length of the edge of the $d$-cube in
	which $\Omega$ is contained. The following dimensionless form of
	the governing equations (\ref{Eq15})-(\ref{Eq17}) (``hat'' is
	dropped for convenience) are available in $\Omega$
	\begin{equation}\label{Eq15a}
		-\nabla\cdot\left(2D(\mathbf{V}_f)+\frac{\lambda_f}{\mu_f}(\nabla\cdot\mathbf{V}_f)\mathbf{I}
		-\phi_fP\mathbf{I}\right)
		+\frac{1}{Da}\mathbf{K}(\boldsymbol{\varsigma})\mathbf{V}_f=\mathbf{b}_f,
	\end{equation}
	\begin{equation}\label{Eq16a}
		-\nabla\cdot\left(\frac{\varrho_t}{(1+\nu_p)}D(\mathbf{U}_s)+\frac{\nu_p\varrho_t}{(1+\nu_p)(1-2\nu_p)}
		(\nabla\cdot\mathbf{U}_s)\mathbf{I}-\phi_sP\mathbf{I}\right)
		-\frac{1}{Da}\mathbf{K}(\boldsymbol{\varsigma})\mathbf{V}_f=\mathbf{b}_s,
	\end{equation}
	\begin{equation}\label{Eq17a}
		\nabla\cdot(\phi_f\mathbf{V}_f)=-\alpha^2_t(1+L_rA_r)(P-1).
	\end{equation}
	In (\ref{Eq15a}) and (\ref{Eq16a}),
	$\mathbf{b}_f$ and $\mathbf{b}_s$ are modified non-dimensional body forces,
	$\varrho_{t}={\mathcal{Y}R^{2}\rho_{f}}/{\mu^{2}_{f}}$ is the dimensionless Young's
	modulus ${Y}$ associated with the solid phase. It contains the response of the solid 
	phase (cellular phase + extracellular matrix) towards viscous drag due to interstitial fluid movement.
	$\alpha^2_{t}=L_P(A/V)\mu_f$ is the strength of solute source, and
	$Da=\frac{K_d\mu_f}{R^2}$ is the Darcy number (permeability
	parameter). The corresponding boundary conditions are
	\begin{equation}\label{Eq18ab}
		\left(2D(\mathbf{V}_f)+\frac{\lambda_f}{\mu_f}
		(\nabla\cdot\mathbf{V}_f)\mathbf{I}-\phi_fP\mathbf{I} \right)\cdot
		\mathbf{n}=\mathbf{T}_\infty~\mbox{and}~\mathbf{U}_s
		=0~~\mbox{on}~~\partial\Omega.
	\end{equation}
	For the sake of writing convenience, set $\lambda=\frac{\lambda_f}{\mu_f},$ $\alpha_1=\frac{\varrho_t}{2(1+\nu_p)},$
	$\alpha_2=\frac{\nu_p\varrho_t}{(1+\nu_p)(1-2\nu_p)},$ and
	$a_0=\alpha^2_t(1+L_rA_r).$ Observe that the system of
	equations (\ref{Eq15a})-(\ref{Eq17a}) is non-linear and fully
	coupled whenever $\boldsymbol{\varsigma}$ is equals to either
	$\mathbf{U}_s$ or $\nabla\cdot\mathbf{U}_s$ which is our primary
	interest for now. The main aim is to study the well-posedness of
	the nonlinear system (\ref{Eq15a})-(\ref{Eq17a}) subject to the
	data (\ref{Eq18ab}).
	\section{Well-posedness of the auxiliary sub-problems}\label{s1}
	$\mathbf{V}_f,$
	$P,$  $\mathbf{U}_s$ are the unknown functions in the system of Eqs (\ref{Eq15a})--(\ref{Eq17a}). We assume the following:
	\begin{description}
		\item[({\bf A})]
		The parameters {$\phi_f>0,$ $\phi_s>0,$ $\lambda\geq0,$
			$\alpha_1>0,$ $\alpha_2>0,$ $a_0>0,$ $Da>0$} are known real
		constants, and the functions $\mathbf{b}_j\in L^2(\Omega)^d$ where
		$j=f,~s,$  $\mathbf{T}_\infty\in L^{2}(\partial\Omega)^d$  are
		also known. $c_k>0,~c_p>0,~c_t>0,~c_s>0$  are some real constants
		that appear in Korn's, Poincare's, trace, and Sobolev's inequalities,
		respectively\footnote{We refer to the Appendix section for details on function spaces and other
			preliminary results.}.
		\item[({\bf B})] Let $\mathbf{K}:\mathbb{R}^d\rightarrow\mathbb{R}^{d^2}$ is a
		symmetric, uniformly bounded, and positive definite matrix. This ensures that
		there exist  positive constants $k_{1}$ and $k_2$ such that for all $\boldsymbol{\xi},~\mathbf{x}\in\mathbb{R}^d$ we have:
		\begin{equation}\label{P1}
			\mbox{(i)}~k_{1}\boldsymbol{\xi}\cdot\boldsymbol{\xi}\leq\mathbf{K}(\mathbf{x})\boldsymbol{\xi}\cdot\boldsymbol{\xi}~~
			\mbox{and}~~\mbox{(ii)}~~||\mathbf{K}(\mathbf{x})||\leq k_2,
		\end{equation}
		\(||\cdot||\) denotes the Euclidean norm.
		\item[({\bf C})] We assume that the hydraulic resistivity 
		$\mathbf{K}$ is Lipschitz continuous. To this extent, let us assume that there
		exists a constant $k_L>0$ such that
		\begin{equation}\label{P3}
			||\mathbf{K}(\mathbf{x})-\mathbf{K}(\mathbf{y})||\leq
			k_L||\mathbf{x}-\mathbf{y}||~\mbox{for all}~\mathbf{x},~\mathbf{y}\in
			\mathbb{R}^d.
		\end{equation}
	\end{description}
	\subsection{Concept of weak formulation}	
	Choose the triplet of test functions \({(\mathbf{W},\mathbf{Z},q)\in H^1(\Omega)^d \times H^1_{0}(\Omega)^d \times L^2(\Omega)}\). {Taking the scalar product of}  \eqref{Eq15a} with \(\mathbf{W}\),  \eqref{Eq16a} with \(\mathbf{Z}\)   and  (\ref{Eq17a}) with \(q,\) and using the boundary conditions \eqref{Eq18ab}, we get the following nonlinear weak formulation: \\   Find the triplet $(\mathbf{V}_f, \mathbf{U}_s,P)\in H^1(\Omega)^d\times H^1_{0}(\Omega)^d\times L^2(\Omega)$ such that
	\begin{align}\nonumber
		2(D(\mathbf{V}_f):D(\mathbf{W}))_{\Omega}+\lambda(\nabla\cdot\mathbf{V}_f,\nabla\cdot\mathbf{W})_{\Omega}
		-\phi_f(P,\nabla\cdot\mathbf{W})_{\Omega}\\\label{WFE1}
		+\frac{1}{Da}(\mathbf{K}(\boldsymbol{\varsigma})\mathbf{V}_f,\mathbf{W})_{\Omega}
		=(\mathbf{b}_f,\mathbf{W})_{\Omega}
		+(\mathbf{T}_\infty,\mathbf{W})_{\partial\Omega}\\ \label{WFE2} 
		\phi_f(\nabla\cdot\mathbf{V}_f,q)_{\Omega}+a_0(P,q)_{\Omega}=(a_0,q)_{\Omega}\\ \nonumber
		2\alpha_1(D(\mathbf{U}_s):D(\mathbf{Z}))_{\Omega}+\alpha_2(\nabla\cdot\mathbf{U}_s,\nabla\cdot\mathbf{Z})_{\Omega}
		-\phi_s(P,\nabla\cdot
		\mathbf{Z})_{\Omega}\\ \label{WFE3}-\frac{1}{Da}(\mathbf{K}(\boldsymbol{\varsigma})\mathbf{V}_f,\mathbf{Z})_{\Omega}=
		(\mathbf{b}_s,\mathbf{Z})_{\Omega}
	\end{align}
	holds for all
	$(\mathbf{W},\mathbf{Z},q)\in H^1(\Omega)^d\times H^1_{0}(\Omega)^d\times L^2(\Omega).$ 
	Here in \(\boldsymbol{\varsigma}\) is   \(\mathbf{U}_s\) and \(\nabla\cdot\mathbf{U}_s.\) 
	\begin{lemma}\label{s2}
		(Equivalence of weak formulations)
		Suppose that parameters and data satisfy assumptions {\bf (A)} and
		{\bf (B)}. Then any solution (in the sense of distributions)
		$(\mathbf{V}_f, \mathbf{U}_s,P)\in H^1(\Omega)^d\times H^1_{0}(\Omega)^d\times L^2(\Omega)$ of the coupled problem
		{(\ref{Eq15a})-(\ref{Eq18ab})} is also a solution to the
		variational problem \eqref{WFE1}-\eqref{WFE3}. Conversely, any solution to the weak
		problem \eqref{WFE1}-\eqref{WFE3} satisfies {(\ref{Eq15a})-(\ref{Eq18ab})}
		in the sense of distributions.
	\end{lemma}
	The proof follows using standard arguments.  We omit to show it.

\subsection{ Decoupled problem corresponding to \eqref{WFE1}-\eqref{WFE3}, Case (a): \(\mathbf{K}(\boldsymbol{\varsigma})=\mathbf{K}(\mathbf{U}_s)\)}
We note that the weak formulation \eqref{WFE1} - \eqref{WFE3} can be decoupled concerning the unknowns $(\mathbf{V}_f, P)$ and $\mathbf{U}_s.$ In case (a) we are dealing with \(\mathbf{K}(\boldsymbol{\varsigma})=\mathbf{K}(\mathbf{U}_s)\) which is a nonlinear function of \(\mathbf{U}_s\)
and satisfies assumptions \eqref{P1} and \eqref{P3}.
In this case, we can
solve \eqref{WFE1} - \eqref{WFE3} sequentially, that is,
given $\boldsymbol{\varsigma}\in H^1_0(\Omega)^d$  find $(\mathbf{V}_f,P)\in H^1(\Omega)^d \times L^2(\Omega)$ such
that
$$(Q_{w_1}(\boldsymbol{\varsigma}))\left\{
\begin{array}{ll}
	2(D(\mathbf{V}_f):D(\mathbf{W}))_{\Omega}+\lambda(\nabla\cdot\mathbf{V}_f,\nabla\cdot\mathbf{W})_{\Omega}
	-\phi_f(P,\nabla\cdot \mathbf{W})_{\Omega}
	\vspace{0.3cm}\\+\frac{1}{Da}(\mathbf{K}(\boldsymbol{\varsigma})\mathbf{V}_f,\mathbf{W})_{\Omega}
	+\phi_f(\nabla\cdot \mathbf{V}_f,q)_{\Omega}
	+a_0(P,q)_{\Omega}\vspace{0.3cm}\\
	=(\mathbf{b}_f,\mathbf{W})_{\Omega}
	+(\mathbf{T}_\infty,\mathbf{W})_{\partial\Omega}+(a_0,q)_{\Omega}
\end{array}
\right.$$ holds for all $(\mathbf{W},q)\in H^1(\Omega)^d \times L^2(\Omega),$
and then for a given pair  $(\mathbf{V}_f,P)\in H^1(\Omega)^d \times L^2(\Omega),$
find $\mathbf{U}_s\in H^1_0(\Omega)^d$ such that
$$(Q_{w_2}(\mathbf{U}_s))\left\{
\begin{array}{ll}
	2\alpha_1(D(\mathbf{U}_s):D(\mathbf{Z}))_{\Omega}+\alpha_2(\nabla\cdot\mathbf{U}_s,\nabla\cdot\mathbf{Z})_{\Omega}
	= \phi_s(P,\nabla\cdot \mathbf{Z})_{\Omega}\vspace{0.3cm}\\+\frac{1}{Da}(\mathbf{K}(\mathbf{U}_s)\mathbf{V}_f,\mathbf{Z})_{\Omega} 
	+(\mathbf{b}_s,\mathbf{Z})_{\Omega} 
\end{array}
\right.$$ holds for all $\mathbf{Z}\in H^1_0(\Omega)^d.$ For notational convenience, we denote combined problem as $(Q_{w_1}(\boldsymbol{\varsigma}))-(Q_{w_2}(\mathbf{U}_s)).$
\subsection{Decoupled problem corresponding to \eqref{WFE1}-\eqref{WFE3}, Case (b): \(\mathbf{K}(\nabla\cdot\boldsymbol{\varsigma})=\mathbf{K}(\nabla\cdot\mathbf{U}_s)\)}
We can
solve \eqref{WFE1} - \eqref{WFE3} sequentially, that is, for a given $\boldsymbol{\varsigma}\in H^1_0(\Omega)^d$  find $(\mathbf{V}_f,P)\in H^1(\Omega)^d \times L^2(\Omega)$ such
that
$$(Q_{w_1}(\nabla\cdot\boldsymbol{\varsigma}))\left\{
\begin{array}{ll}
	2(D(\mathbf{V}_f):D(\mathbf{W}))_{\Omega}+\lambda(\nabla\cdot\mathbf{V}_f,\nabla\cdot\mathbf{W})_{\Omega}
	-\phi_f(P,\nabla\cdot \mathbf{W})_{\Omega}
	\vspace{0.3cm}\\+\frac{1}{Da}(\mathbf{K}(\nabla\cdot\boldsymbol{\varsigma})\mathbf{V}_f,\mathbf{W})_{\Omega}
	+\phi_f(\nabla\cdot \mathbf{V}_f,q)_{\Omega}
	+a_0(P,q)_{\Omega}\vspace{0.3cm}\\
	=(\mathbf{b}_f,\mathbf{W})_{\Omega}
	+(\mathbf{T}_\infty,\mathbf{W})_{\partial\Omega}+(a_0,q)_{\Omega}
\end{array}
\right.$$ holds for all $(\mathbf{W},q)\in H^1(\Omega)^d \times L^2(\Omega),$
and then for a given pair  $(\mathbf{V}_f,P)\in H^1(\Omega)^d \times L^2(\Omega),$
find $\mathbf{U}_s\in H^1_0(\Omega)^d$ such that
$$(Q_{w_2}(\nabla\cdot\mathbf{U}_s))\left\{
\begin{array}{ll}
	2\alpha_1(D(\mathbf{U}_s):D(\mathbf{Z}))_{\Omega}+\alpha_2(\nabla\cdot\mathbf{U}_s,\nabla\cdot\mathbf{Z})_{\Omega}
	= \phi_s(P,\nabla\cdot \mathbf{Z})_{\Omega}\vspace{0.3cm}\\+\frac{1}{Da}(\mathbf{K}(\nabla\cdot\mathbf{U}_s)\mathbf{V}_f,\mathbf{Z})_{\Omega} 
	+(\mathbf{b}_s,\mathbf{Z})_{\Omega} 
\end{array}
\right.$$ holds for all $\mathbf{Z}\in H^1_0(\Omega)^d.$ 
For notational convenience, we denote combined problem as $(Q_{w_1}(\nabla\cdot\boldsymbol{\varsigma}))-(Q_{w_2}(\nabla\cdot\mathbf{U}_s)).$

\subsection{Case (a): Existence and uniqueness results for \((Q_{w_1}(\boldsymbol{\varsigma}))\)}\label{Qa}
In order to solve weak formulation \((Q_{w_1}(\boldsymbol{\varsigma}))\), we use the following method.
We rephrase the weak formulation $(Q_{w_1}(\boldsymbol{\varsigma}))$ into an abstract	setting. Set $\mathbb{Y}=H^1(\Omega)^d \times L^2(\Omega).$ To do so, define a mapping $\mathcal{H}_{\boldsymbol{\varsigma}}$ from 
$\mathbb{Y}$ to $\mathbb{Y}$
by
\begin{align}\nonumber
	\langle\mathcal{H}_{\boldsymbol{\varsigma}}(\mathbf{V}_f,P),(\mathbf{W},q)\rangle=
	2(D(\mathbf{V}_f):D(\mathbf{W}))_{\Omega}+\lambda(\nabla\cdot\mathbf{V}_f,\nabla\cdot\mathbf{W})_{\Omega}
	-\phi_f(P,\nabla\cdot
	\mathbf{W})_{\Omega}\\\label{h1}+\frac{1}{Da}(\mathbf{K}(\boldsymbol{\varsigma})\mathbf{V}_f,\mathbf{W})_{\Omega}
	+\phi_f(\nabla\cdot
	\mathbf{V}_f,q)_{\Omega}+a_0(P,q)_{\Omega}-[(\mathbf{b}_f,\mathbf{W})_{\Omega}
	+(\mathbf{T}_\infty,\mathbf{W})_{\partial\Omega} 
	+(a_0,q)_{\Omega}].
\end{align}
Using the mapping $\mathcal{H}_{\boldsymbol{\varsigma}},$ the variational formulation
\((Q_{w_1}(\boldsymbol{\varsigma}))\) can equivalently be written as: for a given
$\boldsymbol{\varsigma}\in H^1_0(\Omega)^d;$ find $(\mathbf{V}_f,P)\in \mathbb{Y}$ such that
\begin{eqnarray}\label{ef}
	\langle\mathcal{H}_{\boldsymbol{\varsigma}}(\mathbf{V}_f,P),
	(\mathbf{W},q)\rangle=0~~\mbox{for all}~
	(\mathbf{W},q)\in \mathbb{Y}.
\end{eqnarray}
Conversely, if (\ref{ef}) holds then
(\ref{Eq15a}) and (\ref{Eq17a}) with the first boundary condition in the equation 
(\ref{Eq18ab}) satisfy in the sense of distributions. Hence, our
immediate task is to find a pair $(\mathbf{V}_f, P)\in \mathbb{Y}$ that satisfies (\ref{ef}).
In order to do so, we proceed as follows. The mapping $\mathcal{H}_{\boldsymbol{\varsigma}}$ satisfies the following lemma:
\begin{lem}\label{l1} Suppose that parameters and data satisfy assumptions {\bf (A)} and
	{\bf (B)}.
	If ~$\mathcal{H}_{\boldsymbol{\varsigma}}$~ is a mapping from $\mathbb{Y}$ into itself
	defined by (\ref{h1}) then the following statements hold:
	\begin{itemize}
		\item[(i)] $\mathcal{H}_{\boldsymbol{\varsigma}}$ is continuous.
		\item[(ii)] 
		There exists a real number $r>0$ such that
		$$\langle\mathcal{H}_{\boldsymbol{\varsigma}}(\mathbb{V}),\mathbb{V}\rangle>0,~~\mbox{for all}~~\mathbb{V}\in \mathbb{Y}~~
		\mbox{with}~~||\mathbb{V}||_{\mathbb{Y}}=r,$$
	\end{itemize}
	i.e., $\mathcal{H}_{\boldsymbol{\varsigma}}$ is coercive on a ball of radius $r$ in
	$\mathbb{Y}.$ 
	Here, for any
	$\mathbb{V}=(\mathbf{V}_f,P)\in \mathbb{Y}= H^1(\Omega)^d \times L^2(\Omega),$
	$||\cdot||_{\mathbb{Y}}$ is defined as
	$$||\mathbb{V}||^2_{\mathbb{Y}}=||(\mathbf{V}_f,P)||^2_{\mathbb{Y}}
	=||\mathbf{V}_f||^2_{1,\Omega}+||P||^2_{0,\Omega}.$$
\end{lem}

\noindent {\bf Proof :} (i) The continuity of the mapping
$\mathcal{H}_{\boldsymbol{\varsigma}}$ can be shown using the continuity of scalar
product. Indeed, let
$\{\mathbb{V}^m\}_{m\geq1}=
\{(\mathbf{V}^{m}_f,P^{m})\}_{m\geq1}$ be any sequence in $\mathbb{Y}$ that converges strongly to
$\mathbb{V}=(\mathbf{V}_f,P)\in \mathbb{Y}$ as
$m\rightarrow\infty,$ i.e.
\begin{align}\label{c1}
	||\mathbf{V}^{m}_f-\mathbf{V}_f||_{1,\Omega}\rightarrow 0,~~
	||P^{m}-P||_{0,\Omega}\rightarrow 0~\mbox{as }~m\rightarrow\infty.
\end{align}
Relying on the definition of $\mathcal{H}_{\boldsymbol{\varsigma}}$ and on Cauchy-Schwarz
inequality, we get
	\begin{align}\nonumber
		|\langle\mathcal{H}_{\boldsymbol{\varsigma}}(\mathbb{V}^m)
		-\mathcal{H}_{\boldsymbol{\varsigma}}(\mathbb{V}),(\mathbf{W},q)\rangle|\leq
		2||\mathbb{D}^f(\mathbf{V}^m_f-\mathbf{V}_f)||_{0,\Omega}||\mathbb{D}^f(\mathbf{W})||_{0,\Omega} \\\nonumber
		+\lambda ||\nabla\cdot(\mathbf{V}^m_f-\mathbf{V}_f)||_{0,\Omega}||\nabla\cdot\mathbf{W}||_{0,\Omega}
		+\phi_f||P^m-P||_{0,\Omega}||\nabla\cdot\mathbf{W}||_{0,\Omega}\\\nonumber
		+\frac{k_2}{\mbox{Da}}||\mathbf{V}^m_f-\mathbf{V}_f||_{0,\Omega}||\mathbf{W}||_{0,\Omega} 
		+\phi_f||\nabla\cdot(\mathbf{V}^m_f-\mathbf{V}_f)||_{0,\Omega}||q||_{0,\Omega}
		\\\nonumber
		+a_0||P^m-P||_{0,\Omega}||q||_{0,\Omega}.
	\end{align}
Using  
(\ref{c1}), we obtain
\begin{eqnarray}\nonumber
	|\langle\mathcal{H}_{\boldsymbol{\varsigma}}(\mathbf{V}^{m}_f,P^{m})
	-\mathcal{H}_{\boldsymbol{\varsigma}}(\mathbf{V}_f,P)
	,(\mathbf{W},q)\rangle|\rightarrow 0~~ ~\forall~(\mathbf{W},q)~\mbox{as}~
	m\rightarrow\infty.
\end{eqnarray}
This argument establishes the continuity of $\mathcal{H}_{\boldsymbol{\varsigma}}.$\\
(ii) For any
$\mathbb{V}=(\mathbf{V}_f,P)\in \mathbb{Y},$ we have
\begin{align}\nonumber
	\langle\mathcal{H}_{\boldsymbol{\varsigma}}(\mathbb{V}),\mathbb{V}\rangle=
	2||D(\mathbf{V}_f)||^2_{0,\Omega}+\lambda||\nabla\cdot\mathbf{V}_f||^2_{0,\Omega}		-\phi_f(P,\nabla\cdot
	\mathbf{V}_f)_{\Omega}+\frac{1}{Da}(\mathbf{K}(\boldsymbol{\varsigma})\mathbf{V}_f,\mathbf{V}_f)_{\Omega}
	\\\label{h2}
	+\phi_f(\nabla\cdot
	\mathbf{V}_f,P)_{\Omega}+a_0||P||^2_{0,\Omega}
	-[(\mathbf{b}_f,\mathbf{V}_f)_{\Omega}
	+(\mathbf{T}_\infty,\mathbf{V}_f)_{\partial\Omega}+(a_0,P)_{\Omega}].
\end{align}
Using Cauchy-Schwarz, Korn's, and trace inequalities, we
obtain
\begin{align}\nonumber
	\langle\mathcal{H}_{\boldsymbol{\varsigma}}(\mathbb{V}),\mathbb{V}\rangle\geq
	\alpha||\mathbf{V}_f||^2_{1,\Omega}+\lambda||\nabla\cdot\mathbf{V}_f||^2_{0,\Omega}
	+{a_0}||P||^2_{0,\Omega}\\\label{U1}
	-(||\mathbf{b}_f||_{0,\Omega}+\sqrt{c_t}||\mathbf{T}_{\infty}||_{0,\partial\Omega})||\mathbf{V}||_{1,\Omega}
	-||a_0||_{0,\Omega}||P||_{0,\Omega},
\end{align}
where $\alpha=\frac{1}{c_k}\min\{2,\frac{k_1}{Da}\}.$
Further, \eqref{U1} can be rewritten as 
\begin{eqnarray}\label{a3a}
	\langle\mathcal{H}_{\boldsymbol{\varsigma}}(\mathbb{V}),\mathbb{V}\rangle\geq
	\alpha_3(||\mathbf{V}_f||^2_{1,\Omega}
	+||P||^2_{0,\Omega})-\alpha_4||\mathbb{V}||_{1,0,\Omega},
\end{eqnarray}
where
\begin{eqnarray}\label{alp1}
	\alpha_3&=&\min\left\{\alpha,{a_0}\right\},\\\label{alp2}
	\alpha_4&=&[(||\mathbf{b}_f||_{0,\Omega}+\sqrt{c_t}||\mathbf{T}_{\infty}||_{0,\partial\Omega})^2
	+||a_0||^2_{0,\Omega}]^{1/2}.
\end{eqnarray}
If $||\mathbb{V}||_{\mathbb{Y}}=r_0$ for some $r_0>0,$
then we have
\begin{eqnarray}\label{coe}
	\langle\mathcal{H}_{\boldsymbol{\varsigma}}(\mathbb{V}),\mathbb{V}\rangle>0~~~\forall~\mathbb{V}\in \mathbb{Y},
	~~\mbox{when}~~r_0>\frac{\alpha_4}{\alpha_3}.
\end{eqnarray}
This completes the proof of Lemma \ref{l1}.
\begin{rem}
	Note that $\lambda=\frac{\lambda_f}{\mu_f}$ (which is
	the ratio of the two viscosity coefficients) plays a significant role in the coercivity proof. The literature
	\cite{gad1995technical,rajagopal2013new} suggests that researchers
	have debated on the sign (or value) of ${\lambda}$. According to
	the well-known Stokes-hypothesis, $\lambda=-2/3$ \cite{rajagopal2013new}. On the other
	hand, the existing literature also suggests that this ratio can be
	non-negative \cite{gad1995technical}.
	Thus, we consider both of these possibilities.  If $\lambda\geq0$
	then the coercivity of $\mathcal{H}_{\boldsymbol{\varsigma}},$ as shown by us in
	(\ref{coe}), holds. However, when $\lambda<0$ (i.e. a typical
	Stokes hypothesis), the coercivity of $\mathcal{H}_{\boldsymbol{\varsigma}}$ holds with
	relevant restrictions on the constants. For instance,  when $\lambda<0,$ the mapping
	$\mathcal{H}_{\boldsymbol{\varsigma}}$  is coercive if
	$\alpha=\frac{1}{c_k}\min\{2,\frac{k_1}{Da}\}>2/3.$ For
	convenience from here
	onward, we assume $\lambda$ to be a non-negative constant 
\end{rem}
Based on Lemma \ref{l2} (see Appendix \ref{FSpac}), we now present
the following existence results.
\begin{thm}\label{T1} Suppose that {\bf (A)} and {\bf (B)} hold. 
	Then for a given $\boldsymbol{\varsigma}\in H^1_0(\Omega)^d$ the problem (\ref{ef}) has at least one solution
	$(\mathbf{V}_f,P)\in \mathbb{Y}=H^1(\Omega)^d\times L^2(\Omega)$ satisfying the problem
	(\ref{Eq15a}) and (\ref{Eq17a}) with the first boundary condition in the equation 
	(\ref{Eq18ab}) in the sense of distributions. Moreover, the solution $(\mathbf{V}_f,P)$ satisfies the following a priori bound
	\begin{eqnarray}\label{stb}
		||\mathbf{V}_f||^2_{1,\Omega}+||P||^2_{0,\Omega}\leq
		\left(\frac{\alpha_4}{\alpha_3}\right)^{2}.
	\end{eqnarray}
\end{thm}
\noindent {\bf Proof :} To prove this result we use the
Galerkin method. The spaces $H^1(\Omega)^d,$
$L^2(\Omega)$ are separable
Hilbert spaces. Hence, there exist corresponding bases
$\{\mathbf{W}_{i}\}_{i=1}^{\infty}$ and
$\{q_i\}_{i=1}^{\infty}$ of
smooth functions. Let $\mathbb{Y}_m$ be the space spanned by
$\{(\mathbf{W}_{i},q_i)\}_{i=1}^{m}.$ The scalar
product on $\mathbb{Y}_m$ is induced by the scalar product of
$\mathbb{Y}.$ We define the approximate solution
$(\mathbf{V}^m_f,P^m)$ as follows:
\begin{equation}
	\mathbf{V}^m_f=\sum_{i=1}^{m}a_i\mathbf{W}_{i},~~~
	P^m=\sum_{i=1}^{m}c_iq_i,
\end{equation}
with
$$(Q_{m}(\boldsymbol{\varsigma}))\left\{
\begin{array}{ll}\vspace{0.1cm}
	2(D(\mathbf{V}^{m}_f):D(\mathbf{W}))_{\Omega}+\lambda(\nabla\cdot\mathbf{V}^{m}_f,\nabla\cdot\mathbf{W})_{\Omega}
	-\phi_f(P^{m},\nabla\cdot
	\mathbf{W})_{\Omega}\vspace{0.3cm}\\+\frac{1}{Da}(\mathbf{K}(\boldsymbol{\varsigma})\mathbf{V}^{m}_f,\mathbf{W})_{\Omega}
	+\phi_f(\nabla\cdot
	\mathbf{V}^{m}_f,q)_{\Omega}+a_0(P^{m},q)_{\Omega} \vspace{0.3cm}\\=(\mathbf{b}_f,\mathbf{W})_{\Omega}
	+(\mathbf{T}_\infty,\mathbf{W})_{\partial\Omega}+(a_0,q)_{\Omega},
\end{array}
\right.$$ holding for all
$(\mathbf{W},q)\in\mathbb{Y}_m$ with
$a_i,~b_i,~c_i\in\mathbb{R},$ for  $i~=~1,2,\ldots,m.$ The task now is to ensure the existence of
solutions to $(Q_{m}(\boldsymbol{\varsigma}))$ and
show that $(Q_{m}(\boldsymbol{\varsigma}))$ recovers $(Q_{w_1}(\boldsymbol{\varsigma}))$ as
$m\rightarrow\infty.$ The linear structure of $(Q_{m}(\boldsymbol{\varsigma}))$ suggests
that weak convergence is enough to pass the limit. Hence, in
order to do so we define a mapping $\mathcal{H}_m$ inspired by the structure of
mapping $\mathcal{H}$ as
\begin{equation}\label{32}
	\langle\mathcal{H}^m_{\boldsymbol{\varsigma}}(\mathbf{V}_f,P),(\mathbf{W},q)\rangle=
	\langle\mathcal{H}_{\boldsymbol{\varsigma}}(\mathbf{V}_f,P),(\mathbf{W},q)\rangle,~\mbox{for all}~
	(\mathbf{W},q)\in\mathbb{Y}_m
\end{equation}
where $\mathcal{H}_{\boldsymbol{\varsigma}}$ is as defined in (\ref{h1}).  From Lemma
\ref{l1}, we deduce that the mapping $\mathcal{H}^m_{\boldsymbol{\varsigma}}$ satisfies the
conditions needed for Lemma \ref{l2} (see Appendix \ref{FSpac})
and hence, there exists a solution
$(\mathbf{V}^m_f,P^m)\in\mathbb{Y}_m$ for each $m$
such that
\begin{equation}\label{33}
	\langle\mathcal{H}^m_{\boldsymbol{\varsigma}}(\mathbf{V}^m_f,P^m),(\mathbf{W},q)\rangle=0,
	~\mbox{for all}~
	(\mathbf{W},q)\in\mathbb{Y}_m.
\end{equation}
It follows that $(\mathbf{V}^m_f,P^m)$ satisfy $(Q_{m}(\boldsymbol{\varsigma}))$ and $a_i,$ $c_i$ can be determined.\\
{\bf Energy Estimates:} Let $\mathbf{W}=\mathbf{V}^m_f,$
and $q=P^m$ in $(Q_{m}(\boldsymbol{\varsigma}))$  then by
performing calculations similar to those leading to (\ref{a3a}),
we obtain
\begin{eqnarray}\label{a3ab}
	\alpha_3(||\mathbf{V}^m_f||^2_{1,\Omega}
	+||P^m||^2_{0,\Omega})-\alpha_4||(\mathbf{V}^m_f,P^m)||_{1,0,\Omega}\leq0,
\end{eqnarray}
where $\alpha_3$ and $\alpha_4$ are defined in
(\ref{alp1})-(\ref{alp2}). Consequently,
\begin{eqnarray}\label{a3ac}
	||\mathbf{V}^m_f||^2_{1,\Omega}
	+||P^m||^2_{0,\Omega}\leq
	\left(\frac{\alpha_4}{\alpha_3}\right)^{2}.
\end{eqnarray}
Inequality (\ref{a3ac}) implies that the sequence $\{(\mathbf{V}^m_f,P^m)\}_{m\geq1}$ is uniformly bounded in $\mathbb{Y}.$ Hence, it has a subsequence $\{(\mathbf{V}^m_f, P^m)\}_{m\geq1}$ (for convenience, we denote it by the same symbol) and a pair $(\mathbf{V}_f, P)\in\mathbb{Y}$ such that
\begin{eqnarray}\label{a3ad}
	(\mathbf{V}^m_f,P^m)\rightharpoonup(\mathbf{V}_f,P)~~~\mbox{as}~
	m\rightarrow\infty ~~\mbox{weakly~in}~\mathbb{Y}.
\end{eqnarray}
By taking the limit in (\ref{33}) and using the weak convergence
(\ref{a3ad}), we get
\begin{equation}\label{a3e}
	\langle\mathcal{H}_{\boldsymbol{\varsigma}}(\mathbf{V}_f,P),(\mathbf{W},q)\rangle=0,~\mbox{for all}~
	(\mathbf{W},q)\in\mathbb{Y}_m.
\end{equation}
A continuity argument shows that (\ref{a3e}) holds for any $(\mathbf{W},q)\in\mathbb{Y}.$ Hence, $(\mathbf{V}_f, P)$ is a solution of (\ref{ef}) and equivalently, of the weak formulation \((Q_{w_1}(\boldsymbol{\varsigma})).\)  Using the lower semi-continuity property of norm in (\ref{a3ac}), we can achieve the following a priori bound on solution $(\mathbf{V}_f, P)$  given by 
\begin{eqnarray}\label{SB1}
	||\mathbf{V}_f||^2_{1,\Omega}+||P||^2_{0,\Omega}\leq
	\left(\frac{\alpha_4}{\alpha_3}\right)^{2}.
\end{eqnarray}
\begin{prop}\label{T1ba}
	Suppose the hypotheses of Theorem \ref{T1} hold. Then, the weak formulation \((Q_{w_1}(\boldsymbol{\varsigma}))\)  
	has a unique solution that depends continuously on the given data.
\end{prop}
\noindent
{\bf Proof:} {\bf Uniqueness:} Let $(\mathbf{V}^{1}_f,P^{1})$ and
$(\mathbf{V}^{2}_f,P^{2})$ be two solutions that
satisfy equation (\ref{ef}) or equivalently the weak formulation
\((Q_{w_1}(\boldsymbol{\varsigma})).\) Define
$(\mathbf{V}_f,P)=(\mathbf{V}^{1}_f-\mathbf{V}^{2}_f,P^{1}-P^{2}).$
Then from (\ref{ef}), we have
\begin{eqnarray}\label{efg}
	\langle\mathcal{H}_{\boldsymbol{\varsigma}}(\mathbf{V}^{1}_f,P^{1})-\mathcal{H}_{\boldsymbol{\varsigma}}(\mathbf{V}^{2}_f,P^{2}),(\mathbf{W},q)
	\rangle=0
\end{eqnarray}
for all $(\mathbf{W},q)\in\mathbb{Y}.$  Replace
$(\mathbf{W},q)$ by $(\mathbf{V}_f,P)$ in
(\ref{efg}) and using the definition of $\mathcal{H},$ we find
\begin{eqnarray}\nonumber
	\alpha||\mathbf{V}_f||^2_{1,\Omega}+\lambda||\nabla\cdot\mathbf{V}_f||^2_{0,\Omega}
	+{a_0}||P||^2_{0,\Omega}\leq0.
\end{eqnarray}
The above implies
$\mathbf{V}_f=0,$ $\mathbf{U}_s=0$ and $P=0$ a.e. in $\Omega.$
Hence, the weak formulation  \((Q_{w_1}(\boldsymbol{\varsigma}))\)  has a unique solution.\\
{\bf Continuous dependence:}
Let $(\mathbf{V}^1_f,P^1)$ and $(\mathbf{V}^2_f,P^2)$ be two solutions of \((Q_{w_1}(\boldsymbol{\varsigma}))\)  corresponding
to the two sets of {data}
$(\mathbf{T}_{\infty,1},\mathbf{b}_{f,1},a_{0,1})$
and
$(\mathbf{T}_{\infty,2},\mathbf{b}_{f,2},a_{0,2}),$
then the difference
$(\mathbf{V}_f,P)=(\mathbf{V}^{1}_f-\mathbf{V}^{2}_f,P^{1}-P^{2})$
satisfies
\begin{eqnarray}\nonumber
	\alpha||\mathbf{V}_f||^2_{1,\Omega}+\lambda||\nabla\cdot\mathbf{V}_f||^2_{0,\Omega}		+{a_0}||P||^2_{0,\Omega}\\\nonumber
	-(||\mathbf{b}_f||_{0,\Omega}+\sqrt{c_t}||\mathbf{T}_{\infty}||_{0,\partial\Omega})||\mathbf{V}||_{1,\Omega}
	-||a_0||_{0,\Omega}||P||_{0,\Omega}\leq 0,
\end{eqnarray}
or,
\begin{align}\nonumber
	||\mathbf{V}^1_f-\mathbf{V}^2_f||^2_{1,\Omega}
	+||P^1-P^2||^2_{0,\Omega} \\\label{cont}\leq
	\frac{1}{\alpha^2_3}\left[(||\mathbf{b}_{f,1}-\mathbf{b}_{f,2}||_{0,\Omega}
	+\sqrt{c_t}||\mathbf{T}_{\infty,1}-\mathbf{T}_{\infty,2}||_{0,\partial\Omega})^2
	+||a_{0,1}-a_{0,2}||^2_{0,\Omega}\right].
\end{align}
{Thus, if
	$(\mathbf{T}_{\infty,1},\mathbf{b}_{f,1},
	a_{0,1})$ is close to
	$(\mathbf{T}_{\infty,2},\mathbf{b}_{f,2},
	a_{0,2}),$ then the left hand side of (\ref{cont}) (the difference
	of solutions) must be small. This establishes the well-posedness
	of the auxiliary linear problem  \((Q_{w_1}(\boldsymbol{\varsigma})).\) Next we would 
	like to		consider the sub-problem, \((Q_{w_2}(\nabla\cdot\mathbf{U}_s)).\)

	\subsection{Case (a): Existence and uniqueness results for \((Q_{w_2}(\mathbf{U}_s))\)}\label{Qbac}
	
	We note that \(\mathbf{K}(\boldsymbol{\varsigma})=\mathbf{K}(\mathbf{U}_s)\) is a nonlinear function of \(\mathbf{U}_s\) (see assumption (\(A4\)) in subsection \ref{assump}) which makes  \((Q_{w_2}(\mathbf{U}_s))\) a semilinear problem. By introducing a semilinear form \(B(\cdot,\cdot) : H^1_0(\Omega)^d\times H^1_0(\Omega)^d\rightarrow \mathbb{R}\) that is given by 
	\begin{equation}
		B(\mathbf{U}_s,\mathbf{Z})=2\alpha_1(D(\mathbf{U}_s):D(\mathbf{Z}))_{\Omega}+\alpha_2(\nabla\cdot\mathbf{U}_s,\nabla\cdot\mathbf{Z})_{\Omega} -\frac{1}{Da}(\mathbf{K}(\mathbf{U}_s)\mathbf{V}_f,\mathbf{Z})_{\Omega}
	\end{equation}
	and a linear form \(L:H^1_0(\Omega)^d\rightarrow \mathbb{R}\) defined as 
	\begin{equation}
		L(\mathbf{Z})=\phi_s(P,\nabla\cdot \mathbf{Z})_{\Omega}+
		(\mathbf{b}_s,\mathbf{Z})_{\Omega} .
	\end{equation}
	Weak problem \((Q_{w_2}(\mathbf{U}_s))\) can be rewritten as an abstract formulation. For a given pair  $(\mathbf{V}_f,P)\in H^1(\Omega)^d \times L^2(\Omega),$  find $\mathbf{U}_s\in H^1_0(\Omega)^d $ such that  
	\begin{equation}\label{case(b)}
		B(\mathbf{U}_s,\mathbf{Z})=L(\mathbf{Z})~ \mbox{for all}~ \mathbf{Z}\in H^1_0(\Omega)^d.
	\end{equation}
	In order to show existence and uniqueness results for problem \eqref{case(b)}, we will use the Browder-Minty theorem (see Thereom \ref{Minty}, Appendix \ref{FSpac}) which is based on the monotone operator approach. To justify the hypotheses of the Browder-Minty theorem, we prove the following results in the form of lemmas. 
	\begin{lemma}\label{LEM7}
		The correspondence \(\mathbf{Z}\mapsto B(\mathbf{U}_s,\mathbf{Z})\)
		is a bounded linear operator and \({L}\in (H^1_0(\Omega)^d)^*.\)
	\end{lemma}
	{\bf Proof:} Clearly, the mapping  \(\mathbf{Z}\mapsto B(\mathbf{U}_s,\mathbf{Z})\) is linear (obvious) and bounded. Indeed, using Cauchy-Schwarz, H\"{o}lder's, and Sobolev inequalities, we find  
	\begin{align}\nonumber
		|B(\mathbf{U}_s,\mathbf{Z})|\leq 2\alpha_1||D(\mathbf{U}_s)||_{0,\Omega}||D(\mathbf{Z})||_{0,\Omega} +\alpha_2||\nabla\cdot\mathbf{U}_s||_{0,\Omega}||\nabla\cdot\mathbf{Z}||_{0,\Omega} \\ \nonumber+ \frac{1}{Da}||\mathbf{K}(\mathbf{U}_s)||_{L^\infty(\Omega)}||\mathbf{V}_f||_{0,\Omega}||\mathbf{Z}||_{0,\Omega} 
		\leq \Big{(}(2\alpha_1+\alpha_2)||\nabla\mathbf{U}_s||_{0,\Omega} \\ \label{ESTM1} +\frac{k_2\sqrt{c_p}}{Da}||\mathbf{V}_f||_{0,\Omega}\Big{)}||\nabla\mathbf{Z}||_{0,\Omega}.
	\end{align} 
	Now, \(L\) is linear (obvious) and bounded. Indeed, we have 
	\begin{eqnarray}
		|L(\mathbf{Z})|
		\leq (\phi_s||P||_{0,\Omega}+\sqrt{c_p}||\mathbf{b}_s||_{0,\Omega})||\nabla\mathbf{Z}||_{0,\Omega}.
	\end{eqnarray}
	This implies \(L\in (H^1_0(\Omega)^d)^*. \)\\
	The lemma \ref{LEM7} implies that there exists an operator (nonlinear) \(\mathcal{A}: H^1_0(\Omega)^d\rightarrow (H^1_0(\Omega)^d)^*=H^{-1}(\Omega)^d\) with 
	\begin{equation}\label{NON1}
		(\mathcal{A}\mathbf{U}_s,\mathbf{Z})=B(\mathbf{U}_s,\mathbf{Z}).
	\end{equation}
	Thus, the variational problem \eqref{case(b)} equivalently reduces to the operator equation: 
	find \(\mathbf{U}_s\in H^1_0(\Omega)^d \) such that 
	\begin{equation}\label{OE1}
		\mathcal{A}\mathbf{U}_s=L
	\end{equation}
	in the sense 
	\begin{equation}
		(\mathcal{A}\mathbf{U}_s,\mathbf{Z})=L(\mathbf{Z}),~~ \mbox{for all}~\mathbf{Z}\in H^1_0(\Omega)^d.
	\end{equation}
	Further, estimate \eqref{ESTM1} implies the nonlinear operator \(\mathcal{A}\) is bounded. Indeed, 
	\begin{equation}
		||\mathcal{A}\mathbf{U}_s||_{H^{-1}(\Omega)^d} \leq    \left[(2\alpha_1+\alpha_2)||\nabla\mathbf{U}_s||_{0,\Omega}+\frac{k_2\sqrt{c_p}}{Da}||\mathbf{V}_f||_{0,\Omega}\right].
	\end{equation}
	
	\begin{lemma}\label{LEM8}
		If \(\frac{2\alpha_1}{c_k}>\frac{k_Lc_s\sqrt{c_p}\alpha_4}{\alpha_3Da}\) then the semilinear form \(B(\cdot,\cdot)\) is elliptic that is there exists a constant \(c>0\) such that 
		\begin{align}
			B(\mathbf{U}^1_s,\mathbf{U}^1_s-\mathbf{U}^2_s)- B(\mathbf{U}^2_s,\mathbf{U}^1_s-\mathbf{U}^2_s) \geq c||\nabla\mathbf{U}^1_s-\nabla\mathbf{U}^2_s||^2_{0,\Omega}~~\mbox{for all}~\mathbf{U}^1_s,~\mathbf{U}^2_s \in H^1_0(\Omega)^d.
		\end{align}
	\end{lemma}
	Note: This lemma implies \(\mathcal{A}\) is strongly monotone. \\ \vspace{0.3cm}
	{\bf Proof: } Indeed, consider 
	\begin{align}\nonumber
		(\mathcal{A}\mathbf{U}^1_s-\mathcal{A}\mathbf{U}^2_s,
		\mathbf{U}^1_s-\mathbf{U}^2_s)= B(\mathbf{U}^1_s,\mathbf{U}^1_s-\mathbf{U}^2_s)- B(\mathbf{U}^2_s,\mathbf{U}^1_s-\mathbf{U}^2_s)\\ \nonumber
		=  2\alpha_1||D(\mathbf{U}^1_s)-D(\mathbf{U}^2_s)||^2_{0,\Omega}+ 
		\alpha_2||\nabla\cdot\mathbf{U}^1_s-\nabla\cdot\mathbf{U}^2_s||^2_{0,\Omega}
		\\ \nonumber -\frac{1}{Da}((\mathbf{K}(\mathbf{U}^1_s)-\mathbf{K}(\mathbf{U}^2_s))\mathbf{V}_f,\mathbf{U}^1_s-\mathbf{U}^2_s)_{\Omega}
		\geq 2\alpha_1||D(\mathbf{U}^1_s)-D(\mathbf{U}^2_s)||^2_{0,\Omega} \\ \nonumber+ 
		\alpha_2||\nabla\cdot\mathbf{U}^1_s-\nabla\cdot\mathbf{U}^2_s||^2_{0,\Omega}
		-\frac{1}{Da}||\mathbf{K}(\mathbf{U}^1_s)-\mathbf{K}(\mathbf{U}^2_s)||_{0,\Omega}||\mathbf{V}_f||_{L^4(\Omega)}||\mathbf{U}^1_s-\mathbf{U}^2_s||_{L^4(\Omega)} \\ \nonumber
		\geq 
		\frac{2\alpha_1}{c_k}||\nabla\mathbf{U}^1_s-\nabla\mathbf{U}^2_s||^2_{0,\Omega} 
		-\frac{k_Lc_s\sqrt{c_p}\alpha_4}{\alpha_3Da}||\nabla\mathbf{U}^1_s-\nabla\mathbf{U}^2_s||^2_{0,\Omega} \\ \label{Monotone} 
		=
		\left(\frac{2\alpha_1}{c_k}
		-\frac{k_Lc_s\sqrt{c_p}\alpha_4}{\alpha_3Da}\right)||\nabla\mathbf{U}^1_s-\nabla\mathbf{U}^2_s||^2_{0,\Omega}
	\end{align}
	To reach \eqref{Monotone}, we have applied H\"{o}lder's, Poincare's and Sobolev's inequalities and Lipschitz continuous property of \(\mathbf{K}.\) Thus, if \(\frac{2\alpha_1}{c_k}>\frac{k_Lc_s\sqrt{c_p}\alpha_4}{\alpha_3Da}\) then \(B\) is elliptic with \(c=\left(\frac{2\alpha_1}{c_k}
	-\frac{k_Lc_s\sqrt{c_p}\alpha_4}{\alpha_3Da}\right).\)
	\begin{lemma}\label{LEM8b}
		The nonlinear operator \(\mathcal{A}\) as in \eqref{NON1} is continuous from \(H^1_0(\Omega)^d\) to \(H^{-1}(\Omega)^d.\)
	\end{lemma}
	{\bf Proof:} Let 
	\(\mathbf{U}^n_s \rightarrow \mathbf{U}_s\) in \(H^1_0(\Omega)^d\) as \(n\rightarrow \infty.\) 
	Consider
	\begin{align}\nonumber
		|(\mathcal{A}\mathbf{U}^n_s-\mathcal{A}\mathbf{U}_s,\mathbf{Z})| \leq 
		2\alpha_1||D(\mathbf{U}^n_s)-D(\mathbf{U}_s)||_{0,\Omega}||D(\mathbf{Z})||_{0,\Omega} \\ \nonumber+ 
		\alpha_2||\nabla\cdot\mathbf{U}^n_s-\nabla\cdot\mathbf{U}_s||_{0,\Omega}
		||\nabla\cdot\mathbf{Z}||_{0,\Omega}  +\frac{1}{Da}||\mathbf{K}(\mathbf{U}^n_s)-\mathbf{K}(\mathbf{U}_s)||_{0,\Omega}||\mathbf{V}_f||_{L^4(\Omega)}||\mathbf{Z}||_{L^4(\Omega)}
	\end{align}
	On applying Lipschitz continuity of \(\mathbf{K}\) and Sobolev inequality, we find
	\begin{align}\nonumber
		||\mathcal{A}\mathbf{U}^n_s-\mathcal{A}\mathbf{U}_s||_{H^{-1}(\Omega)^d} \leq 
		2\alpha_1||D(\mathbf{U}^n_s)-D(\mathbf{U}_s)||_{0,\Omega}+ 
		\alpha_2||\nabla\cdot\mathbf{U}^n_s-\nabla\cdot\mathbf{U}_s||_{0,\Omega}
		\\ \nonumber +\frac{k_Lc_s}{Da}||\mathbf{U}^n_s-\mathbf{U}_s||_{0,\Omega}
	\end{align}
	and  \(||\mathcal{A}\mathbf{U}^n_s-\mathcal{A}\mathbf{U}_s||_{H^{-1}(\Omega)^d}\rightarrow 0\) as \(n\rightarrow \infty.\) That is, \(\mathcal{A}\) is continuous. 
	
	\begin{thm}\label{T1ab} Suppose that assumptions {\bf (A)},  {\bf (B)} and {\bf (C)} hold. 
		Further, if the given data and non-dimensional parameters satisfy the following assumption 
		\begin{equation}\label{INEQ2}
			\frac{2\alpha_1}{c_k}>\frac{k_Lc_s\sqrt{c_p}\alpha_4}{\alpha_3Da},
		\end{equation}
		then for a given pair $(\mathbf{V}_f,P)\in H^1(\Omega)^d \times L^2(\Omega)^d$
		the problem \((Q_{w_2}(\mathbf{U}_s))\) has a unique solution 
		\(\mathbf{U}_s\in H^1_0(\Omega)^d. \) 
		Moreover, if 
		\begin{equation}\label{INEQ3}
			\frac{2\alpha_1}{c_k}>\frac{k_2\sqrt{c_p}\alpha_4}{\alpha_3Da}
		\end{equation}
		holds then
		\(\mathbf{U}_s\) satisfies the following a priori estimate
		\begin{eqnarray}\label{stbb1}
			||\nabla\mathbf{U}_s||_{0,\Omega} \leq \frac{1}{\left(\frac{2\alpha_1}{c_k}-\frac{k_2\sqrt{c_p}\alpha_4}{\alpha_3Da}\right)}
			\left(\frac{\phi_s\alpha_4}{\alpha_3}+\sqrt{c_p}||\mathbf{b}_s||_{0,\Omega}\right).
		\end{eqnarray} 
	\end{thm}
	{\bf Proof:}
	The analysis shown in Lemmas \ref{LEM8} and \ref{LEM8b} justify that \(\mathcal{A}\) satisfies the hypothesis of the Browder-Minty theorem (see Theorem \ref{Minty}, Appendix \ref{FSpac}).
	Consequently, the operator equation \eqref{OE1} or the problem \eqref{case(b)} has a unique solution \(\mathbf{U}_s\in H^1_0(\Omega)^d\) for any given pair $(\mathbf{V}_f,P)\in H^1(\Omega)^d \times L^2(\Omega).$  Further, if \eqref{INEQ3} holds
	then \(\mathbf{U}_s\) satisfies the apriory estimate \eqref{stbb}. Indeed, from  \eqref{case(b)} replacing \(\mathbf{Z}=\mathbf{U}_s\)  and using Poincare's inequality, we have
	\begin{eqnarray}\nonumber
		B(\mathbf{U}_s,\mathbf{U}_s)=L(\mathbf{U}_s)\leq (\phi_s||P||_{0,\Omega}+\sqrt{c_p}||\mathbf{b}_s||_{0,\Omega})||\nabla\mathbf{U}_s||_{0,\Omega}
	\end{eqnarray}
	Making use of the definition of \(B(\mathbf{U}_s,\mathbf{U}_s)\) and Korn's, H\"{o}lder's and Sobolev inequalities and boundedness property of   \(\mathbf{K}\), we obtain  
	\begin{eqnarray}\nonumber
		||\nabla\mathbf{U}_s||_{0,\Omega} \leq \frac{1}{\left(\frac{2\alpha_1}{c_k}-\frac{k_2\sqrt{c_p}\alpha_4}{\alpha_3Da}\right)}
		\left(\frac{\phi_s\alpha_4}{\alpha_3}+\sqrt{c_p}||\mathbf{b}_s||_{0,\Omega}\right).
	\end{eqnarray}
	The analysis in subsections \ref{Qa}-\ref{Qbac} describes the existence and uniqueness of problems $(Q_{w_1}(\boldsymbol{\varsigma}))-(Q_{w_2}(\mathbf{U}_s)).$  In the next subsections, we focus on developing existence and uniqueness results corresponding to $(Q_{w_1}(\nabla\cdot\boldsymbol{\varsigma}))-(Q_{w_2}(\nabla\cdot\mathbf{U}_s))$ that is the case (b).


	\subsection{Case (b): Existence and uniqueness of a solution to \((Q_{w_1}(\nabla\cdot\boldsymbol{\varsigma}))\)}\label{Qba}
	Analysis in this subsection is analogous to the arguments in the subsection \ref{Qa}. Thus, we only state the main theorem and omit the proof. 
	
	\begin{thm}\label{T10} Suppose that the assumptions {\bf (A)} and {\bf (B)} hold. 
		Then for a given $\boldsymbol{\varsigma}\in H^1_0(\Omega)^d$ the problem \((Q_{w_1}(\nabla\cdot\boldsymbol{\varsigma}))\) has a unique solution
		$(\mathbf{V}_f,P)\in \mathbb{Y}=H^1(\Omega)^d\times L^2(\Omega)$ that depends continuously on the given data. Moreover, the solution $(\mathbf{V}_f,P)$ satisfies the following a priori bound
		\begin{eqnarray}\label{stb0}
			||\mathbf{V}_f||^2_{1,\Omega}+||P||^2_{0,\Omega}\leq
			\left(\frac{\alpha_4}{\alpha_3}\right)^{2}.
		\end{eqnarray}
	\end{thm}
	{\bf Proof:} Proof of this theorem follows from Theorem \ref{T1} and Proposition \ref{T1ba}.
	\subsection{Case (b): Existence and uniqueness of a solution to \((Q_{w_2}(\nabla\cdot\mathbf{U}_s))\)}\label{Qbb}
		
We note that \(\mathbf{K}(\nabla\cdot\mathbf{U}_s)\) is a nonlinear function of \(\nabla\cdot\mathbf{U}_s\) (see assumption (\(A4\)) in subsection \ref{assump}) which makes  \((Q_{w_2}(\nabla\cdot\mathbf{U}_s))\) a semilinear problem. By introducing a semilinear form \(B(\cdot,\cdot) : H^1_0(\Omega)^d\times H^1_0(\Omega)^d\rightarrow \mathbb{R}\) that is given by 
\begin{equation}
	B(\mathbf{U}_s,\mathbf{Z})=2\alpha_1(D(\mathbf{U}_s):D(\mathbf{Z}))_{\Omega}+\alpha_2(\nabla\cdot\mathbf{U}_s,\nabla\cdot\mathbf{Z})_{\Omega} -\frac{1}{Da}(\mathbf{K}(\nabla\cdot\mathbf{U}_s)\mathbf{V}_f,\mathbf{Z})_{\Omega}
\end{equation}
and a linear form \(L:H^1_0(\Omega)^d\rightarrow \mathbb{R}\) defined as 
\begin{equation}
	L(\mathbf{Z})=\phi_s(P,\nabla\cdot \mathbf{Z})_{\Omega}+
	(\mathbf{b}_s,\mathbf{Z})_{\Omega} .
\end{equation}
Weak problem \((Q_{w_2}(\nabla\cdot\mathbf{U}_s))\) can be rewritten as an abstract formulation. For a given pair  $(\mathbf{V}_f,P)\in H^1(\Omega)^d \times L^2(\Omega),$  find $\mathbf{U}_s\in H^1_0(\Omega)^d $ such that  
\begin{equation}\label{case(b)9}
	B(\mathbf{U}_s,\mathbf{Z})=L(\mathbf{Z})~ \mbox{for all}~ \mathbf{Z}\in H^1_0(\Omega)^d.
\end{equation}
Similar to subsection \ref{Qbac} to show existence and uniqueness results for problem \eqref{case(b)9}, we will use the Browder-Minty theorem (see Theorem \ref{Minty}, Appendix \ref{FSpac}) which is based on the monotone operator approach. We state the main theorem without proof which can be done similar to the proof of Theorem \ref{T1ab}.

\begin{thm}\label{T1ab9} Suppose that assumptions {\bf (A)},  {\bf (B)} and {\bf (C)} hold. 
	Further, if the given data and non-dimensional parameters satisfy the following assumption 
	\begin{equation}\label{INEQ29}
		\frac{2\alpha_1}{c_k}>\frac{k_Lc_s\alpha_4}{\alpha_3Da},
	\end{equation}
	then for a given pair $(\mathbf{V}_f,P)\in H^1(\Omega)^d \times L^2(\Omega)^d$
	the problem \((Q_{w_2}(\nabla\cdot\mathbf{U}_s))\) has a unique solution 
	\(\mathbf{U}_s\in H^1_0(\Omega)^d. \) 
	Moreover, if 
	\begin{equation}\label{INEQ39}
		\frac{2\alpha_1}{c_k}>\frac{k_2\sqrt{c_p}\alpha_4}{\alpha_3Da}
	\end{equation}
	holds then
	\(\mathbf{U}_s\) satisfies the following a priori estimate
	\begin{eqnarray}\label{stbb19}
		||\nabla\mathbf{U}_s||_{0,\Omega} \leq \frac{1}{\left(\frac{2\alpha_1}{c_k}-\frac{k_2\sqrt{c_p}\alpha_4}{\alpha_3Da}\right)}
		\left(\frac{\phi_s\alpha_4}{\alpha_3}+\sqrt{c_p}||\mathbf{b}_s||_{0,\Omega}\right).
	\end{eqnarray} 
\end{thm}
 {\bf Proof:} The proof of this theorem can be done similarly to the proof of Theorem \ref{T1ab}.
	The analysis in subsections \ref{Qba}-\ref{Qbb} completes the existence and uniqueness of solution to problem $(Q_{w_1}(\nabla\cdot\boldsymbol{\varsigma}))-(Q_{w_2}(\nabla\cdot\mathbf{U}_s)).$  In the next section, we focus on the development of existence and uniqueness results corresponding to coupled nonlinear problems  $(Q_{w_1}(\mathbf{U}_s))-(Q_{w_2}(\mathbf{U}_s))$ and $(Q_{w_1}(\nabla\cdot\mathbf{U}_s))-(Q_{w_2}(\nabla\cdot\mathbf{U}_s)),$  respectively, by converting them into a fixed-point problem. 
	
	\section{Case (a): Reduction to a fixed-point problem for $(Q_{w_1}(\mathbf{U}_s))-(Q_{w_2}(\mathbf{U}_s))$ } \label{case(abb)} 
	
	We note that for a given \(\boldsymbol{\varsigma}\in H^1_0(\Omega)^d,\) the problem \((Q_{w_1}(\boldsymbol{\varsigma}))\) has a unique solution \((\mathbf{V}_f, P)\in H^1(\Omega)^d\times L^2(\Omega)\) (see subsection \ref{Qa}). Consequently, we can define a mapping \(T_1:  H^1_0(\Omega)^d\rightarrow H^1(\Omega)^d\times L^2(\Omega)\) such that \(T_1(\boldsymbol{\varsigma})=(\mathbf{V}_f, P).\) Further, for a given pair \((\mathbf{V}_f, P)\in H^1(\Omega)^d\times L^2(\Omega),\) the problem \((Q_{w_2}(\mathbf{U}_s))\) has a unique solution \(\mathbf{U}_s\in H^1_0(\Omega)^d\) (see subsection \ref{Qbac}). Therefore, we can define a mapping \(T_2: H^1(\Omega)^d\times L^2(\Omega)\rightarrow H^1_0(\Omega)^d\) such that \(T_2(\mathbf{V}_f, P)=\mathbf{U}_s.\) Now, in order to get the fixed-point problem corresponding to $(Q_{w_1}(\mathbf{U}_s))-(Q_{w_2}(\mathbf{U}_s)),$ we define a composition map \(T=T_2\circ T_1: H^1_0(\Omega)^d \rightarrow H^1_0(\Omega)^d\) such that 
	\begin{equation}\label{FPM1}
		T(\boldsymbol{\varsigma})=(T_2\circ T_1)(\boldsymbol{\varsigma})=T_2(T_1(\boldsymbol{\varsigma}))=T_2(\mathbf{V}_f, P)=\mathbf{U}_s.
	\end{equation}
	Thus, a fixed-point of mapping \(T\) solves the coupled nonlinear problem 
	\((Q_{w_1}(\mathbf{U}_s))-(Q_{w_2}(\mathbf{U}_s))\) or equivalently, \eqref{WFE1} - \eqref{WFE3} when \(\boldsymbol{\varsigma}=\mathbf{U}_s.\) In order to show that the mapping \(T\) has a fixed point, we use Schauder's fixed-point theorem (see Thereom \ref{Schauder}, Appendix \ref{FSpac})\footnote{For the fixed point approach, we are inspired by the working techniques used in \cite{camano2016augmented,caucao2020fully} }. Thus, in the following analysis, we prove some results in the form of lemmas to verify the hypotheses of Schauder's fixed-point theorem.
	
	\subsection{Analysis of the fixed point problem}\label{AFP2}
	Throughout this subsection, we assume that hypotheses of Theorem \ref{T1}, Proposition \ref{T1ba} and Theorem \ref{T1ab} hold. Then, we have
	\begin{lemma}\label{LEM5}
		Given \(r>0,\) let \(\mathbf{W}_r\) be a closed and convex subset of \(H^1(\Omega)^d\) defined by 
		\[\mathbf{W}_r=\{\mathbf{w}\in H^1_0(\Omega)^d : ||\mathbf{w}||_{H^1_0(\Omega)^d}=||\nabla\mathbf{w}||_{0,\Omega}\leq r\},\] and assume that the data satisfies
		\begin{equation}\label{F1}
			\frac{1}{\left(\frac{2\alpha_1}{c_k}- \frac{k_2\sqrt{c_p}\alpha_4}{\alpha_3Da}\right)}\left(\phi_s\frac{\alpha_4}{\alpha_3}+\sqrt{c_p}||\mathbf{b}_s||_{0,\Omega}\right)\leq r.  
		\end{equation}
		Then \({T}(\mathbf{W}_r)\subseteq \mathbf{W}_r.\)
	\end{lemma}
	{Proof: } For any \(\mathbf{w}\in \mathbf{W}_r,\) using estimate \eqref{stbb1}, we find
	\begin{align}\nonumber
		||T(\mathbf{w})||_{H^1_0(\Omega)^d}=||T_2(T_1(\mathbf{w}))||_{H^1_0(\Omega)^d}=||T_2(\mathbf{V}_f,P)||_{H^1_0(\Omega)^d}=||\mathbf{U}_s||_{H^1_0(\Omega)^d}\\ \label{F2}=
		||\nabla\mathbf{U}_s||_{0,\Omega}\leq \frac{1}{\left(\frac{2\alpha_1}{c_k}- \frac{k_2\sqrt{c_p}\alpha_4}{\alpha_3Da}\right)}\left(\phi_s\frac{\alpha_4}{\alpha_3}+\sqrt{c_p}||\mathbf{b}_s||_{0,\Omega}\right).
	\end{align}
	The above estimate \eqref{F2} together with the assumption \eqref{F1} imply \(T(\mathbf{w}) \in \mathbf{W}_r,\) which proves \(T(\mathbf{W}_r)\subseteq \mathbf{W}_r.\)
	\vspace{0.3cm}
	\begin{lemma}\label{LEM6}
		The map \({T}_1 : H^1_0(\Omega)^d \rightarrow H^1(\Omega)^d\times L^2(\Omega)\) satisfies 
		\begin{eqnarray}\label{F1ab}
			||{T}_1(\boldsymbol{\varsigma})-{T}_1(\tilde{\boldsymbol{\varsigma}})||_{\mathbb{Y}} \leq \frac{\sqrt{2}k_L\alpha_4c_s}{\alpha^2_3Da}||\boldsymbol{\varsigma}
			-\tilde{\boldsymbol{\varsigma}}||_{0,\Omega}.
		\end{eqnarray}
	\end{lemma}
	{\bf Proof:}
	Given \(\boldsymbol{\varsigma}, \tilde{\boldsymbol{\varsigma}}\in H^1_0(\Omega)^d ,\)  let 
	\((\mathbf{V}_f, P),~(\tilde{\mathbf{V}}_f, \tilde{P})\in H^1(\Omega)^d\times L^2(\Omega)\) be the corresponding solutions of \((Q_{w_1}(\boldsymbol{\varsigma})).\) Then, equation \eqref{ef} implies
	\begin{eqnarray}
		\langle\mathcal{H}_{\boldsymbol{\varsigma}}(\mathbf{V}_f,P),(\mathbf{W},q)\rangle=0~\mbox{for all}~
		(\mathbf{W},q)\in\mathbb{Y},\\ 
		\langle\mathcal{H}_{\tilde{\boldsymbol{\varsigma}}}(\tilde{\mathbf{V}}_f, \tilde{P}),(\mathbf{W},q)\rangle=0~\mbox{for all}~
		(\mathbf{W},q)\in\mathbb{Y}.
	\end{eqnarray}
	Taking the difference between the above equations, we get 
	\begin{eqnarray}
		\langle\mathcal{H}_{\boldsymbol{\varsigma}}(\mathbf{V}_f,P)-\mathcal{H}_{\tilde{\boldsymbol{\varsigma}}}(\tilde{\mathbf{V}}_f, \tilde{P}),(\mathbf{W},q)\rangle=0~\mbox{for all}~
		(\mathbf{W},q)\in\mathbb{Y}.
	\end{eqnarray}
	Replacing \(\mathbf{W}=\mathbf{V}_f-\tilde{\mathbf{V}}_f\) and \(q=P-\tilde{P}\) in the above equation and using the definition of \(\mathcal{H}_{\boldsymbol{\varsigma}}\) and \(\mathcal{H}_{\tilde{\boldsymbol{\varsigma}}},\) we get
	\begin{eqnarray}\nonumber
		2||D(\mathbf{V}_f)-D(\tilde{\mathbf{V}}_f)||^2_{0,\Omega}+\lambda||\nabla\cdot\mathbf{V}_f-\nabla\cdot\tilde{\mathbf{V}}_f||^2_{0,\Omega}
		+a_0||P-\tilde{P}||^2_{0,\Omega}\\=-\frac{1}{Da}(\mathbf{K}(\boldsymbol{\varsigma})(\mathbf{V}_f-\tilde{\mathbf{V}}_f),\mathbf{V}_f-\tilde{\mathbf{V}}_f)_{\Omega}-\frac{1}{Da}((\mathbf{K}(\boldsymbol{\varsigma})-\mathbf{K}(\tilde{\boldsymbol{\varsigma}}))\tilde{\mathbf{V}}_f,\mathbf{V}_f-\tilde{\mathbf{V}}_f)_{\Omega}.
	\end{eqnarray}
	Using Korn's, H\"{o}lder's inequalities and  estimate \eqref{SB1}, we obtain
	\begin{eqnarray}\nonumber
		||\mathbf{V}_f-\tilde{\mathbf{V}}_f||^2_{1,\Omega}+||P-\tilde{P}||^2_{0,\Omega}
		\leq \frac{\alpha_4}{\alpha^2_3 Da}||\mathbf{K}(\boldsymbol{\varsigma}) 
		-\mathbf{K}(\tilde{\boldsymbol{\varsigma}})||_{0,\Omega}||\mathbf{V}_f-\tilde{\mathbf{V}}_f||_{L^4(\Omega)}.
	\end{eqnarray}
	Further, using Lipschitz continuous property of \(\mathbf{K}\) and Sobolev's inequality, we get
	\begin{eqnarray}\nonumber
		\left(||\mathbf{V}_f-\tilde{\mathbf{V}}_f||^2_{1,\Omega}+||P-\tilde{P}||^2_{0,\Omega}\right)
		\leq \frac{k_L\alpha_4{c_s}}{\alpha^2_3Da}||\boldsymbol{\varsigma}
		-\tilde{\boldsymbol{\varsigma}}||_{0,\Omega}||\mathbf{V}_f-\tilde{\mathbf{V}}_f||_{1,\Omega},
	\end{eqnarray}
	or,
	\begin{align}\label{F4bb}
		||(\mathbf{V}_f, P)-(\tilde{\mathbf{V}}_f, \tilde{P})||_{\mathbb{Y}}=\left(||\mathbf{V}_f-\tilde{\mathbf{V}}_f||^2_{1,\Omega}+||P-\tilde{P}||^2_{0,\Omega}\right)^{1/2}
		\leq \frac{\sqrt{2}k_L\alpha_4c_s}{\alpha^2_3Da}||\boldsymbol{\varsigma}
		-\tilde{\boldsymbol{\varsigma}}||_{0,\Omega}.
	\end{align}
	This establishes \eqref{F1ab}.
	
	\begin{lemma}\label{LEM9} 
		The map \({T}_2 : H^1(\Omega)^d\times L^2(\Omega)\rightarrow H^1_0(\Omega)^d\) satisfies
		\begin{align}
			||{T}_2(\mathbf{V}_f, P)-{T}_2(\tilde{\mathbf{V}}_f, \tilde{P})||_{H^1_0(\Omega)}  \leq {\beta}\left[||P-\tilde{P}||_{0,\Omega}
			+||\mathbf{V}_f-\tilde{\mathbf{V}}_f||_{1,\Omega}\right]
		\end{align}  
		where \({\beta}= \frac{1}{\left(\frac{2\alpha_1}{c_k}-\frac{k_Lc_s\alpha_4\sqrt{c_p}}{\alpha_3Da}\right)}\max\left\{\phi_s, \frac{k_2\sqrt{c_p}}{Da}\right\}.\)
	\end{lemma}
	{\bf Proof:} 
	Given \((\mathbf{V}_f, P),~(\tilde{\mathbf{V}}_f, \tilde{P})\in H^1(\Omega)^d\times L^2(\Omega),\)  let \(\mathbf{U}_s,~\tilde{\mathbf{U}}_s\in H^1_0(\Omega)^d  \) be the corresponding solutions of \((Q_{w_2}(\mathbf{U}_s)).\) Then, \eqref{case(b)} implies
	\begin{eqnarray}
		B(\mathbf{U}_s,\mathbf{Z})=L(\mathbf{Z})~ \mbox{for all}~ \mathbf{Z}\in H^1_0(\Omega)^d,\\
		B(\tilde{\mathbf{U}}_s,\mathbf{Z})=L(\mathbf{Z})~ \mbox{for all}~ \mathbf{Z}\in H^1_0(\Omega)^d.
	\end{eqnarray}
	Taking the difference between the above equations, we get 
	\begin{align}\nonumber
		2\alpha_1(D(\mathbf{U}_s)-D(\tilde{\mathbf{U}}_s):D(\mathbf{Z}))_{\Omega}+\alpha_2(\nabla\cdot\mathbf{U}_s-\nabla\cdot\tilde{\mathbf{U}}_s,\nabla\cdot\mathbf{Z})_{\Omega}\\ \nonumber
		=\phi_s(P-\tilde{P},\nabla\cdot \mathbf{Z})_{\Omega}
		+\frac{1}{Da}(\mathbf{K}(\mathbf{U}_s)(\mathbf{V}_f-\tilde{\mathbf{V}}_f),\mathbf{Z})_{\Omega}
		\\ \nonumber+\frac{1}{Da}((\mathbf{K}(\mathbf{U}_s)-\mathbf{K}(\tilde{\mathbf{U}}_s))\tilde{\mathbf{V}}_f,\mathbf{Z})_{\Omega}.
	\end{align}
	Replace \(\mathbf{Z}=\mathbf{U}_s -\tilde{\mathbf{U}}_s\) and using Cauchy-Schwarz and H\"{o}lder's inequalities, we get
	\begin{align}\nonumber
		2\alpha_1||D(\mathbf{U}_s)-D(\tilde{\mathbf{U}}_s)||^2_{0,\Omega}+\alpha_2||\nabla\cdot\mathbf{U}_s-\nabla\cdot\tilde{\mathbf{U}}_s||^2_{0,\Omega}
		\leq \phi_s||P-\tilde{P}||_{0,\Omega}||\nabla\cdot (\mathbf{U}_s -\tilde{\mathbf{U}}_s)||_{0,\Omega} \\ \nonumber
		+\frac{1}{Da}||\mathbf{K}(\mathbf{U}_s)||_{L^\infty(\Omega)}||\mathbf{V}_f-\tilde{\mathbf{V}}_f||_{0,\Omega}||\mathbf{U}_s -\tilde{\mathbf{U}}_s||_{0,\Omega} \\ \nonumber
		+\frac{1}{Da}||\mathbf{K}(\mathbf{U}_s)-\mathbf{K}(\tilde{\mathbf{U}}_s)||_{0,\Omega}||\tilde{\mathbf{V}}_f||_{L^4(\Omega)}||\mathbf{U}_s -\tilde{\mathbf{U}}_s||_{L^4(\Omega)}.
	\end{align}
	Moreover, using Lipschitz property of \(\mathbf{K}\), Korn's, Poincare's, and Sobolev's embedding inequalities and estimate \eqref{SB1}, we get 
	\begin{align}\nonumber
		||T_2(\mathbf{V}_f,P)-T_2(\tilde{\mathbf{V}}_f,\tilde{P})||_{H^1_0(\Omega)^d}=||\mathbf{U}_s-\tilde{\mathbf{U}}_s||_{H^1_0(\Omega)^d}=
		||\nabla\mathbf{U}_s -\nabla\tilde{\mathbf{U}}_s||_{0,\Omega}\\ \nonumber
		\leq \frac{1}{\left(\frac{2\alpha_1}{c_k}-\frac{k_Lc_s\alpha_4\sqrt{c_p}}{\alpha_3Da}\right)}\left[\phi_s||P-\tilde{P}||_{0,\Omega}
		+\frac{k_2\sqrt{c_p}}{Da}||\mathbf{V}_f-\tilde{\mathbf{V}}_f||_{0,\Omega}\right],
	\end{align}
	or,   
	\begin{align}\label{F5ab}
		||T_2(\mathbf{V}_f,P)-T_2(\tilde{\mathbf{V}}_f,\tilde{P})||_{H^1_0(\Omega)^d}
		\leq \beta\left[||P-\tilde{P}||_{0,\Omega}
		+||\mathbf{V}_f-\tilde{\mathbf{V}}_f||_{1,\Omega}\right]
	\end{align}
	where \({\beta}= \frac{1}{\left(\frac{2\alpha_1}{c_k}-\frac{k_Lc_s\alpha_4\sqrt{c_p}}{\alpha_3Da}\right)}\max\left\{\phi_s, \frac{k_2\sqrt{c_p}}{Da}\right\}.\)
	\begin{lemma}\label{LEM10}
		The map  \({T}={T}_2\circ{T}_1 : H^1_0(\Omega)^d \subset L^2(\Omega)^d\rightarrow H^1_0(\Omega)^d\) satisfies
		\begin{align}\label{F1bb}
			||{T}(\boldsymbol{\varsigma})-{T}(\tilde{\boldsymbol{\varsigma}})||_{H^1_0(\Omega)^d}\leq 
			\frac{2{\beta} k_L\alpha_4{c_s}}{\alpha^2_3Da}||\boldsymbol{\varsigma}
			-\tilde{\boldsymbol{\varsigma}}||_{0,\Omega}\leq \frac{2{\beta} k_L\alpha_4c_s\sqrt{c_p}}{\alpha^2_3Da}||\nabla\boldsymbol{\varsigma}
			-\nabla\tilde{\boldsymbol{\varsigma}}||_{0,\Omega}.
		\end{align}
	\end{lemma}
	{\bf Proof:} 
	From \eqref{F5ab}, we have 
	\begin{align}\nonumber
		||T(\boldsymbol{\varsigma})-T(\tilde{\boldsymbol{\varsigma}})||_{H^1_0(\Omega)^d)}= ||T_2(\mathbf{V}_f,P)-T_2(\tilde{\mathbf{V}}_f,\tilde{P})||_{H^1_0(\Omega)^d}
		\\\label{F6ab} \leq \beta\left[||P-\tilde{P}||_{0,\Omega}
		+||\mathbf{V}_f-\tilde{\mathbf{V}}_f||_{1,\Omega}\right]
	\end{align}
	We achieve \eqref{F1bb} with the help of \eqref{F4bb} and \eqref{F6ab}.
	
	\begin{thm}\label{THM2}
		The mapping  \({T} : \mathbf{W}_r\subset H^1_0(\Omega)^d \rightarrow \mathbf{W}_r\subset H^1_0(\Omega)^d \) is continuous and
		\(\overline{{T}(\mathbf{W}_r)}\) is compact.
	\end{thm}
	{\bf Proof:} The continuity of \(T\) follows in a straightforward manner from \eqref{F1bb}.
	Now given a sequence \(\{\boldsymbol{\varsigma}_k\}_{k\in\mathbb{N}}\) of \(\mathbf{W}_r\) which is clearly bounded, there exists a subsequence \(\{\boldsymbol{\varsigma}_{k_j}\}\subset \{\boldsymbol{\varsigma}_k\}_{k\in\mathbb{N}} \) and \(\boldsymbol{\varsigma}\in H^1_0(\Omega)^d\) such that \(\boldsymbol{\varsigma}_{k_j}\rightharpoonup \boldsymbol{\varsigma}\) in \(H^1_0(\Omega)^d.\) In this way, thanks to the compact embedding of \(H^1_0(\Omega)^d\) in \(L^2(\Omega)^d,\) which implies \(\boldsymbol{\varsigma}_{k_j} \rightarrow \boldsymbol{\varsigma}\) in \(L^2(\Omega)^d,\) which combined with \eqref{F1bb} implies that \(T(\boldsymbol{\varsigma}_{k_j})\rightarrow T(\boldsymbol{\varsigma}).\) This proves that \(\overline{T(\mathbf{W}_r)}\) is compact.
	
	\begin{thm} \label{THM4}
		Suppose that the hypotheses of Theorem \ref{T1}, Proposition \ref{T1ba} and Theorem \ref{T1ab} hold. Then, the mapping \(T\) defined in \ref{FPM1} has a fixed point \(\mathbf{U}_s\in H^1_0(\Omega)^d\) which in turn implies that coupled problem \((Q_{w_1}(\mathbf{U}_s))-(Q_{w_2}(\mathbf{U}_s))\) has a solution \((\mathbf{V}_f, P)\in H^1(\Omega)^d\times L^2(\Omega)\) and \( \mathbf{U}_s\in H^1_0(\Omega)^d.\) Further, if 
		\begin{equation}\label{INEQ4}
			\frac{2{\beta} k_L\alpha_4{c_s}\sqrt{c_p}}{\alpha^2_3Da}<1
		\end{equation}
		then \(T\) has a unique fixed point \( \mathbf{U}_s\in H^1_0(\Omega)^d\) which in turn implies that coupled problem \((Q_{w_1}(\mathbf{U}_s))-(Q_{w_2}(\mathbf{U}_s))\) has a unique solution \((\mathbf{V}_f, P)\in H^1(\Omega)^d\times L^2(\Omega)\) and \( \mathbf{U}_s\in H^1_0(\Omega)^d.\)
	\end{thm}
	{\bf Proof:} 
	The above Theorem \ref{THM2} implies that \(T\) satisfies all hypotheses of Schauder's fixed-point theorem (see Thereom \ref{Schauder}, Appendix \ref{FSpac}). Consequently, the mapping \(T\) has a fixed-point \(\mathbf{U_s}\in H^1_0(\Omega).\) This means that coupled problem \((Q_{w_1}(\mathbf{U}_s))-(Q_{w_2}(\mathbf{U}_s))\) or equivalently, \eqref{WFE1} - \eqref{WFE3} when \(\boldsymbol{\varsigma}=\mathbf{U}_s\) has a solution \((\mathbf{V}_f, P)\in H^1(\Omega)^d\times L^2(\Omega)\) and \( \mathbf{U}_s\in H^1_0(\Omega)^d.\) Further, if \ref{INEQ4} is true then \(T: \mathbf{W}_r \rightarrow \mathbf{W}_r \) is a strict contraction mapping. That implies \(T\) has a unique fixed point \( \mathbf{U}_s\in H^1_0(\Omega)^d\) due to Banach's fixed-point theorem (see p.415 \cite{salsa2016partial}). This means that coupled problem \((Q_{w_1}(\mathbf{U}_s))-(Q_{w_2}(\mathbf{U}_s))\) or equivalently, \eqref{WFE1} - \eqref{WFE3} when \(\boldsymbol{\varsigma}=\mathbf{U}_s\) has a unique solution \((\mathbf{V}_f, P)\in H^1(\Omega)^d\times L^2(\Omega)\) and \( \mathbf{U}_s\in H^1_0(\Omega)^d.\)
	\vspace{0.3cm}
	
		
	\begin{rem}	 
			{({\bf Continuous dependence})}	If the non-dimensional parameters and constants satisfy the following assumption
			\begin{align}\label{P4}
				2\alpha^* Da>\frac{k_L\alpha_4c_s}{\alpha_3},~
				\frac{4\alpha_1Da}{c_k}>\left(\frac{c_pk^2_2}{k_1}+\frac{k_L\alpha_4c_s(c_p+{2}\sqrt{c_p})}{\alpha_3}\right),~
				2\alpha_2\geq\frac{\phi^2_s}{a_0},
			\end{align}
			where $\alpha^*=\frac{1}{c_k}\min\{2,\frac{k_1}{2Da}\}.$
			Then the continuous dependence (\ref{cont}) holds. Indeed, we have
			\begin{align}\nonumber
				||\mathbf{V}^1_f-\mathbf{V}^2_f||^2_{1,\Omega}
				+||\nabla(\mathbf{U}^1_s- \mathbf{U}^2_s)||^2_{0,\Omega}
				+||P^1-P^2||^2_{0,\Omega} \leq
				\frac{1}{\alpha^2_6}[(||\mathbf{b}_{f,1}-\mathbf{b}_{f,2}||_{0,\Omega}
				\\\label{cont1}
				+\sqrt{c_t}||\mathbf{T}_{\infty,1}-\mathbf{T}_{\infty,2}||_{0,\partial\Omega})^2
				+||a_{0,1}-a_{0,2}||^2_{0,\Omega}+c_p||\mathbf{b}_{s,1}-\mathbf{b}_{s,2}||^2_{0,\Omega}],
			\end{align}
			where
			$\alpha_5=\min\left\{\alpha^*-\frac{k_L\alpha_4c_s}{2\alpha_3Da},
			\left(\frac{2\alpha_1}{c_k}-\frac{c_pk^2_2}{2k_1Da}-\frac{k_L\alpha_4c_s(c_p+{2}\sqrt{c_p})}{2\alpha_3Da}
			\right), \frac{a_0}{2}\right\}.$ 
		\end{rem}

		\section{Case (b): Reduction to a fixed-point problem to $(Q_{w_1}(\nabla\cdot\mathbf{U}_s))-(Q_{w_2}(\nabla\cdot\mathbf{U}_s))$  }\label{case(ab)}
		We note that for a given \(\boldsymbol{\varsigma}\in H^1_0(\Omega)^d,\) the problem \((Q_{w_1}(\nabla\cdot\boldsymbol{\varsigma}))\) has a unique solution \((\mathbf{V}_f, P)\in H^1(\Omega)^d\times L^2(\Omega)\) (see subsection \ref{Qba}). Consequently, we can define a mapping \(T_1: H^1_0(\Omega)^d \rightarrow H^1(\Omega)^d\times L^2(\Omega)\) such that \(T_1(\boldsymbol{\varsigma})=(\mathbf{V}_f, P).\) Further, for a given pair \((\mathbf{V}_f, P)\in H^1(\Omega)^d\times L^2(\Omega),\) the problem \((Q_{w_2}(\nabla\cdot\mathbf{U}_s))\) has a unique solution \(\mathbf{U}_s\in H^1_0(\Omega)^d\) (see subsection \ref{Qbb}). Therefore, we can define a mapping \(T_2: H^1(\Omega)^d\times L^2(\Omega)\rightarrow H^1_0(\Omega)^d\) such that \(T_2(\mathbf{V}_f, P)=\mathbf{U}_s.\) Now, in order to get the fixed-point problem corresponding to $(Q_{w_1}(\nabla\cdot\mathbf{U}_s))-(Q_{w_2}(\nabla\cdot\mathbf{U}_s)),$ we define a composition map \(T=T_2\circ T_1: H^1_0(\Omega)^d \rightarrow H^1_0(\Omega)^d\) such that 
		\begin{equation}\label{FPM} 
			T(\boldsymbol{\varsigma})=T_2(T_1(\boldsymbol{\varsigma}))=T_2(\mathbf{V}_f, P)=\mathbf{U}_s.
		\end{equation}
		Thus, a fixed-point of mapping \(T\) solves the coupled nonlinear problem 
		\((Q_{w_1}((\nabla\cdot\mathbf{U}_s))\) and \((Q_{w_2}((\nabla\cdot\mathbf{U}_s))\) or equivalently, \eqref{WFE1} - \eqref{WFE3} when \(\boldsymbol{\varsigma}=\nabla\cdot\mathbf{U}_s.\) In order to show that the mapping \(T\) has a fixed point, we use Schauder's fixed-point theorem (see Thereom \ref{Schauder}, Appendix \ref{FSpac}).
		Thus, in the following analysis, we state some results in the form of lemmas to verify the hypotheses of Schauder's fixed-point theorem. The proof of the following results is almost similar to the proof presented in subsection \ref{AFP2}, so we omit the details.
		\subsection{Analysis of the fixed point problem}\label{AFP1}
		Throughout this subsection, we assume that hypotheses of Theorems \ref{T10} and \ref{T1ab9} hold. Then, we have
		\begin{lemma}\label{LEM1}
			Given \(r>0,\)  let \(\mathbf{W}_r\) be a closed and convex subset of \(H^1(\Omega)^d\) defined by 
			\[\mathbf{W}_r=\{\mathbf{w}\in H^1_0(\Omega)^d : ||\mathbf{w}||_{H^1_0(\Omega)^d}=||\nabla\mathbf{w}||_{0,\Omega}\leq r\}\] and assume that the data satisfy
			\begin{equation}\label{F1aa}
				\frac{1}{\left(\frac{2\alpha_1}{c_k}-\frac{k_2\sqrt{c_p}\alpha_4}{\alpha_3Da}\right)}
		\left(\frac{\phi_s\alpha_4}{\alpha_3}+\sqrt{c_p}||\mathbf{b}_s||_{0,\Omega}\right)\leq r.  
			\end{equation}
			Then \(T(\mathbf{W}_r)\subseteq \mathbf{W}_r.\)
		\end{lemma}
  {\bf Proof: } This lemma is proved similarly to lemma  \ref{LEM5}.
		\vspace{0.3cm}
		\begin{lemma}\label{LEM2}
			The map \(T_1 : H^1_0(\Omega)^d \rightarrow H^1(\Omega)^d\times L^2(\Omega)\) satisfies 
			\begin{eqnarray}\label{F1a}
				||T_1(\boldsymbol{\varsigma})-T_1(\tilde{\boldsymbol{\varsigma}})||_{\mathbb{Y}} \leq \frac{\sqrt{2}k_L\alpha_4{c_s}}{\alpha^2_3Da}||\nabla\cdot\boldsymbol{\varsigma}
				-\nabla\cdot\tilde{\boldsymbol{\varsigma}}||_{0,\Omega},
			\end{eqnarray}
			(see Lemma \ref{l1} for the definition of the space \(\mathbb{Y}\)).
		\end{lemma}
		{\bf Proof: } This lemma is proved similarly to lemma  \ref{LEM6}.
		\begin{lemma}\label{LEM3}
			The map \(T_2 : H^1(\Omega)^d\times L^2(\Omega)\rightarrow H^1_0(\Omega)^d\) satisfies
			\begin{align}
				||T_2(\mathbf{V}_f, P)-T_2(\tilde{\mathbf{V}}_f, \tilde{P})||_{H^1_0(\Omega)^d}  \leq \tilde{\beta}\left[||P-\tilde{P}||_{0,\Omega}
				+||\mathbf{V}_f-\tilde{\mathbf{V}}_f||_{1,\Omega}\right],
			\end{align}  
			where  \(\tilde{\beta}= \frac{1}{\left(\frac{2\alpha_1}{c_k}-\frac{k_L\alpha_4c_s}{\alpha_3Da}\right)}\max\left\{\phi, \frac{k_2\sqrt{c_p}}{Da}
			\right\}\).
		\end{lemma}
  {\bf Proof: } This lemma is proved similarly to lemma  \ref{LEM9}.

		\begin{lemma}\label{LEM4}
			The map  \(T=T_2\circ T_1 : H^1_0(\Omega)^d \rightarrow H^1_0(\Omega)^d\) satisfies
			\begin{align}\label{F5a}
				||T(\boldsymbol{\varsigma})-T(\tilde{\boldsymbol{\varsigma}})||_{H^1_0(\Omega)^d)}
				\leq  \frac{2\tilde{\beta} k_L\alpha_4{c_s}}{\alpha^2_3Da}||\nabla\cdot\boldsymbol{\varsigma}
				-\nabla\cdot\tilde{\boldsymbol{\varsigma}}||_{0,\Omega}\leq  \frac{2 \tilde{\beta} k_L\alpha_4c_s}{\alpha^2_3Da}||\nabla\boldsymbol{\varsigma}
				-\nabla\tilde{\boldsymbol{\varsigma}}||_{0,\Omega}.
			\end{align}
		\end{lemma}
		{\bf Proof: } This lemma is proved similarly to lemma  \ref{LEM10}.
		
		\begin{thm}\label{THM6} 
			Suppose that the hypotheses of Theorems \ref{T10} and \ref{T1ab9} hold and 
   if the following assumption
			\begin{equation}\label{INEQ5}
				\frac{2 \tilde{\beta} k_L\alpha_4c_s}{\alpha^2_3Da}<1
			\end{equation}
			holds then \(T\)  has a unique fixed point \( \mathbf{U}_s\in H^1_0(\Omega)^d\) which in turn implies that coupled problem \((Q_{w_1}(\nabla\cdot\mathbf{U}_s))-(Q_{w_2}(\nabla\cdot\mathbf{U}_s))\) has a unique solution \((\mathbf{V}_f, P)\in H^1(\Omega)^d\times L^2(\Omega)\) and \( \mathbf{U}_s\in H^1_0(\Omega)^d.\)
		\end{thm}
		{\bf Proof:}
  If \( \frac{4 \tilde{\beta} k_L\alpha_4c_s}{\alpha^2_3Da}<1\) then lemma \ref{LEM4} implies \(T: \mathbf{W}_r\subset H^1_0(\Omega)^d \rightarrow \mathbf{W}_r\subset H^1_0(\Omega)^d \) is a strict contraction mapping. That implies \(T\) has a unique fixed point \( \mathbf{U}_s\in H^1_0(\Omega)^d\) due to Banach's fixed-point theorem. That means the coupled problem $(Q_{w_1}(\nabla\cdot\mathbf{U}_s))-(Q_{w_2}(\nabla\cdot\mathbf{U}_s))$  or equivalently, the problem \eqref{WFE1} - \eqref{WFE3} with \(\boldsymbol{\varsigma}=\nabla\cdot\mathbf{U}_s\) has a unique solution \((\mathbf{V}_f, P)\in H^1(\Omega)^d\times L^2(\Omega)\) and \( \mathbf{U}_s\in H^1_0(\Omega)^d.\) 
		
		\begin{rem}\noindent
			\begin{itemize}
				\item If the non-dimensional parameters and the constants satisfy the following constraints
				\begin{align}\label{P14}
					2\alpha^* Da>\frac{k_Lc_s\alpha_4}{\alpha_3},
					\frac{4\alpha_1Da}{c_k}>\left(\frac{c_pk^2_2}{k_1}
					+\frac{3k_Lc_s\alpha_4}{\alpha_3}\right),
					2\alpha_2\geq\frac{\phi^2_s}{a_0}.
				\end{align} then the continuous dependence (\ref{cont1}) holds with the modified constant
				$$\alpha_6=\min\left\{\alpha^*-\frac{k_Lc_s\alpha_4}
				{2\alpha_3Da},\left(\frac{2\alpha_1}{c_k} -\frac{c_pk^2_2}{2k_1Da} -\frac{3k_Lc_s\alpha_4}{2\alpha_3Da}
				\right), \frac{a_0}{2}\right\}.$$
			\end{itemize}
		\end{rem}

		\section{Unbounded $\mathbf{K}$}\label{Unbound}
		One may note that Theorems in sections \ref{s1}, \ref{case(abb)} and \ref{case(ab)} are proven under the
		boundedness assumption  (ii) of {(\ref{P1})}
		and Lipschitz continuity {(\ref{P3})}
		of $\mathbf{K}=\mathbf{K}(\boldsymbol{\varsigma}).$ 
		In this section, we would like to relax such assumptions for case (a)  \(\boldsymbol{\varsigma}=\mathbf{U}_s\) and case (b) \(\boldsymbol{\varsigma}=\nabla\cdot\mathbf{U}_s\).
		Instead of boundedness property (ii) of (\ref{P1}) we assume there exists a constant $k_0>0,$ such that following
		sub-linear growth condition holds
		\begin{eqnarray}\label{P5}
			{||\mathbf{K}(\mathbf{x})||\leq
				k_0(1+||\mathbf{x}||)~\mbox{for all}~\mathbf{x}\in\mathbb{R}^d}.
		\end{eqnarray}
  The above growth condition yields for case (a)
		\begin{eqnarray}
   \label{GC1}
			||\mathbf{K}(\mathbf{U}_s)||_{0,\Omega}\leq
			\sqrt{2}k_0(\sqrt{|\Omega|}+||\mathbf{U}_s||_{0,\Omega}).
		\end{eqnarray}
		We have the following result for case (a).
		\begin{thm}\label{THM3}
			Assume that the data and parameters satisfy the assumption
			{\bf (A)}, $(i)$ of (\ref{P1}) and (\ref{P3}), \eqref{INEQ2}, \eqref{INEQ3} and \(\mathbf{K}\) satisfies the growth condition \eqref{GC1}. Further, if the following parameter constraint
			\begin{align}\label{Rassup1}
				\frac{2\alpha_1}{c_k}>\frac{\sqrt{2c_p}\,c_sk_0\alpha_4}{\alpha_3Da}		
			\end{align}
			is satisfied, then  \eqref{WFE1} - \eqref{WFE3} together with \(\boldsymbol{\varsigma}=\mathbf{U}_s\)
			has a
			solution $(\mathbf{V}_f, \mathbf{U}_s,P)\in H^1(\Omega)^d \times H^1_0(\Omega)^d \times L^2(\Omega)$
			such that
			\begin{align}
				||(\mathbf{V}_f,P)||_{\mathbb{Y}} \leq\frac{\alpha_4}{\alpha_3}
			\end{align}
			and
			\begin{align}
				||\nabla\mathbf{U}_s||_{0,\Omega}\leq \frac{1}{{\left(\frac{2\alpha_1}{c_k}-\frac{\sqrt{2c_p}\,c_sk_0\alpha_4}{\alpha_3Da}\right)}}\left[
				\sqrt{c_p}||\mathbf{b}_s||_{0,\Omega}+\frac{\alpha_4}{\alpha_3}\left(\phi_s+\frac{\sqrt{2|\Omega|}\,c_sk_0}{Da}\right)\right]
			\end{align}
			where $|\Omega|$ denotes the area or volume of the domain $\Omega.$
		\end{thm} 
		\noindent {\bf Proof:} We are inspired by the working techniques from \cite{ccecsmeliouglu2012existence}. Here, we are dealing with particular nonlinear structures of the hydraulic
		resistivity. We approximate the operators $\mathbf{K}$ with a sequence of uniformly positive definite, bounded operators,
		$\{\mathbf{K}^m\}_{m\geq1}$ defined by
		\begin{eqnarray}\label{P6}
			\mathbf{K}^m_{ij}(\mathbf{U}_s):
			=\min(m,\mathbf{K}_{ij}(\mathbf{U}_s)), ~\mbox{for all}~m\in\mathbb{N}.
		\end{eqnarray}
		{Since $\mathbf{K}^m$ is bounded for each $m,$ Theorem
			\mbox{\ref{THM4}} implies that there exists a triplet\\ $(\mathbf{V}^m_f,
			\mathbf{U}^m_s,P^m)$ in $H^1(\Omega)^d \times H^1_0(\Omega)^d \times L^2(\Omega)$
			such that}
			\begin{align}\nonumber
				2(D(\mathbf{V}^{m}_f):D(\mathbf{W}))_{\Omega}+\lambda(\nabla\cdot\mathbf{V}^{m}_f,\nabla\cdot\mathbf{W})_{\Omega}
				-\phi_f(P^{m},\nabla\cdot \mathbf{W})_{\Omega}
				\\ \nonumber+\frac{1}{Da}(\mathbf{K}^m(\mathbf{U}^{m}_s)\mathbf{V}^{m}_f,\mathbf{W})_{\Omega}
				+\phi_f(\nabla\cdot \mathbf{V}^{m}_f,q)_{\Omega}
				+a_0(P^{m},q)_{\Omega}\\ \label{QUAAA}
				=(\mathbf{b}_f,\mathbf{W})_{\Omega}
				+(\mathbf{T}_\infty,\mathbf{W})_{\partial\Omega}+(a_0,q)_{\Omega},
			\end{align}
			\begin{align}\nonumber
				2\alpha_1(D(\mathbf{U}^{m}_s):D(\mathbf{Z}))_{\Omega}
				+\alpha_2(\nabla\cdot\mathbf{U}^{m}_s,\nabla\cdot\mathbf{Z})_{\Omega}
				= \phi_s(P^{m},\nabla\cdot \mathbf{Z})_{\Omega}
				\\ \label{QUBBB}+\frac{1}{Da}(\mathbf{K}^m(\mathbf{U}^{m}_s)\mathbf{V}^{m}_f,\mathbf{Z})_{\Omega}+
				(\mathbf{b}_s,\mathbf{Z})_{\Omega}.
			\end{align}
		We	desire a uniform bound for the sequence $\{(\mathbf{V}^m_f,
		\mathbf{U}^m_s,P^m)\}$ which is independent of $m.$ In order to do
		so, replace $(\mathbf{W},q)$ by $(\mathbf{V}^m_f,P^m)$ in \eqref{QUAAA}. Using the positive
		definiteness of $\mathbf{K}^m$ and the Cauchy-Schwarz, trace, Korn's inequalities,
		we get the following estimate
		\begin{align}\label{P9}
			||(\mathbf{V}^m_f,P^m)||_{\mathbb{Y}} \leq\frac{\alpha_4}{\alpha_3}.
		\end{align}
		One can observe the right hand side of (\ref{P9}) is independent of $m.$
		Next, replacing $\mathbf{Z}$ by $\mathbf{U}^m_s$ in \eqref{QUBBB} and use H\"{o}lder's, Poincare's, Korn's, Sobolev inequalities and   \eqref{P5},  to get 
		\begin{align}\nonumber
			\frac{2\alpha_1}{c_k}||\nabla\mathbf{U}^m_s||^2_{0,\Omega}
			\leq \sqrt{c_p}||\mathbf{b}_s||_{0,\Omega}||\nabla\mathbf{U}^m_s||_{0,\Omega}+\phi_s||P^m||_{0,\Omega}||\nabla\cdot\mathbf{U}^m_s||_{0,\Omega}
			\\ \nonumber+\frac{1}{Da}||\mathbf{K}(\mathbf{U}^{m}_s)||_{0,\Omega}||\mathbf{V}^{m}_f||_{L^4(\Omega)}||\mathbf{U}^{m}_s||_{L^4(\Omega)}
		\end{align}
		\begin{align}\nonumber
			\frac{2\alpha_1}{c_k}||\nabla\mathbf{U}^m_s||^2_{0,\Omega}
			\leq \sqrt{c_p}||\mathbf{b}_s||_{0,\Omega}||\nabla\mathbf{U}^m_s||_{0,\Omega}+\phi_s||P^m||_{0,\Omega}||\nabla\mathbf{U}^m_s||_{0,\Omega}
			\\ \nonumber +\frac{c_s}{Da}[\sqrt{2}k_0(\sqrt{\Omega}+\sqrt{c_p}||\nabla\mathbf{U}^m_s||_{0,\Omega})]||\mathbf{V}^{m}_f||_{1,\Omega}||\nabla\mathbf{U}^{m}_s||_{0,\Omega}
		\end{align}
		\begin{align}\label{P10}
			||\nabla\mathbf{U}^m_s||_{0,\Omega}\leq \frac{1}{{\left(\frac{2\alpha_1}{c_k}-\frac{\sqrt{2c_p}\,c_sk_0\alpha_4}{\alpha_3Da}\right)}}\left[
			\sqrt{c_p}||\mathbf{b}_s||_{0,\Omega}+\frac{\alpha_4}{\alpha_3}\left(\phi_s+\frac{\sqrt{2|\Omega|}\,c_sk_0}{Da}\right)\right]
		\end{align}
		The estimates
		(\ref{P9}) and (\ref{P10}) imply that $(\mathbf{V}^m_f,
		\mathbf{U}^m_s,P^m)\in H^1(\Omega)^d \times H^1_0(\Omega)^d \times L^2(\Omega)$
		uniformly bounded for all $m\geq1.$ Hence, there exists a triplet
		$(\mathbf{V}_f, \mathbf{U}_s,P)\in H^1(\Omega)^d \times H^1_0(\Omega)^d \times L^2(\Omega)$ and a sub-sequence of $(\mathbf{V}^m_f,
		\mathbf{U}^m_s, P^m)$ (we denote by the same symbol) such that
		\begin{eqnarray}\label{P11}
			(\mathbf{V}^m_f, \mathbf{U}^m_s,P^m)\rightharpoonup(\mathbf{V}_f,
			\mathbf{U}_s,P)~\mbox{weakly in}~H^1(\Omega)^d \times H^1_0(\Omega)^d \times L^2(\Omega),
		\end{eqnarray}
		and the compact embedding $H^1(\Omega)^d\hookrightarrow
		L^4(\Omega)^d$ yields
		\begin{eqnarray}\label{P12}
			(\mathbf{V}^m_f, \mathbf{U}^m_s)\rightarrow(\mathbf{V}_f,
			\mathbf{U}_s)~\mbox{strongly in}~L^4(\Omega)^d\times
			L^4(\Omega)^d.
		\end{eqnarray}
		We seek to pass the limit in
		\eqref{QUAAA}-\eqref{QUBBB} as $m\rightarrow\infty.$ We observe that
		the weak convergence (\ref{P11}) is sufficient to pass the limit
		in the linear terms of \eqref{QUAAA}-\eqref{QUBBB} however the
		nonlinear terms demand strong convergence (\ref{P12}).
		Subsequently, the nonlinear terms of the problem
		\eqref{QUAAA}-\eqref{QUBBB} which can be rewritten as \\
		(i)  for \eqref{QUAAA}
		\begin{align}\label{chap2cw4ba}
			\underbrace{((\mathbf{K}^m(\mathbf{U}^{m}_s)-\mathbf{K}(\mathbf{U}_s))\mathbf{V}^{m}_f,\mathbf{W})_{\Omega}}
			+\underbrace{(\mathbf{K}(\mathbf{U}_s)(\mathbf{V}^{m}_f-\mathbf{V}_f),\mathbf{W})_{\Omega}}+
			(\mathbf{K}(\mathbf{U}_s)\mathbf{V}_f,\mathbf{W})_{\Omega}
		\end{align}
		(ii) for \eqref{QUBBB}
		\begin{align}\label{chap2cw4ab}
			\underbrace{((\mathbf{K}^m(\mathbf{U}^m_s)-\mathbf{K}(\mathbf{U}_s))\mathbf{V}^{m}_f,\mathbf{Z})_{\Omega}}
			+\underbrace{(\mathbf{K}(\mathbf{U}_s)(\mathbf{V}^{m}_f-\mathbf{V}_f),\mathbf{Z})_{\Omega}}
			+(\mathbf{K}(\mathbf{U}_s)\mathbf{V}_f, \mathbf{Z})_{\Omega},
		\end{align}
		{Since $\mathbf{U}^m_s$ converges to $\mathbf{U}_s$ strongly in
			$L^2(\Omega)^d,$ it implies $\mathbf{U}^{m}_s-\mathbf{U}_s\rightarrow0$ a.e. in $\Omega$ up to a
			sub-sequence.  The fact that $\mathbf{K}$ is Lipschitz guarantees that} 
		\begin{eqnarray}\label{cw3n}
			{\mathbf{K}^m(\mathbf{U}^{m}_s)-\mathbf{K}(\mathbf{U}_s)\rightarrow0~\mbox{a.e.
					in}~~\Omega.}
		\end{eqnarray}
		{Indeed, we have}
		\begin{align}\nonumber
			{|\mathbf{K}^m(\mathbf{U}^{m}_s)(\mathbf{x})-\mathbf{K}(\mathbf{U}_s)(\mathbf{x})|}	{=|\mathbf{K}^m(\mathbf{U}^{m}_s)(\mathbf{x})-\mathbf{K}^m(\mathbf{U}_s)(\mathbf{x})+	\mathbf{K}^m(\mathbf{U}_s)(\mathbf{x})+\mathbf{K}(\mathbf{U}_s)(\mathbf{x})|}\\\nonumber				
			{\leq k_L|\mathbf{U}^{m}_s(\mathbf{x})-\mathbf{U}_s(\mathbf{x})|
				+|\mathbf{K}^m(\mathbf{U}_s)(\mathbf{x})-\mathbf{K}(\mathbf{U}_s)(\mathbf{x})|.}
		\end{align}
		{As $\mathbf{U}^{m}_s-\mathbf{U}_s\rightarrow0$ a.e. in
			$\Omega$ together with the definition of the truncation function
			implies				$\mathbf{K}^m(\mathbf{U}_s)-\mathbf{K}(\mathbf{U}_s)\rightarrow0$
			a.e. in $\Omega.$ This, establishes \mbox{(\ref{cw3n})}.}
		Using the bound (\ref{P9}) and the convergence result (\ref{cw3n}), we get
		\begin{align}
			\lim_{m\rightarrow\infty}((\mathbf{K}^m(\mathbf{U}^{m}_s)-\mathbf{K}(\mathbf{U}_s))
			\mathbf{V}^{m}_f,\mathbf{W})_{\Omega}=0,~~\lim_{m\rightarrow\infty}((\mathbf{K}^m(\mathbf{U}^{m}_s)-\mathbf{K}(\mathbf{U}_s))
			\mathbf{V}^{m}_f,\mathbf{Z})_{\Omega}=0.
		\end{align}		
		The bound (\ref{P10})
  together with estimate \eqref{GC1}, and the strong convergence $\mathbf{V}^m_f\rightarrow\mathbf{V}_f$ in
		$L^4(\Omega)^d$ leads to 
		\begin{align}	
			\lim_{m\rightarrow\infty}(\mathbf{K}(\mathbf{U}_s)
			(\mathbf{V}^{m}_f-\mathbf{V}_f),\mathbf{W})_{\Omega}=0, ~~   \lim_{m\rightarrow\infty}(\mathbf{K}(\mathbf{U}_s)
			(\mathbf{V}^{m}_f-\mathbf{V}_f),\mathbf{Z})_{\Omega}=0
		\end{align}	
		Thus, the terms with under braces in \eqref{chap2cw4ba} and \eqref{chap2cw4ab} tend to zero as $m\rightarrow\infty.$ Hence, we can say	that \eqref{QUAAA}-\eqref{QUBBB} recovers \eqref{WFE1}-\eqref{WFE3} as  $m\rightarrow\infty.$ This completes the proof of the Theorem \ref{THM3}.\\ 

		\begin{rem}

  \begin{itemize}
      \item
			We note that in the above case i.e. when \(\mathbf{K}(\mathbf{U}_s)\) is not bounded, however, satisfies the sub-linear growth condition \eqref{P5}, the uniqueness of solutions  holds
			whenever the non-dimensional parameters and constants satisfy the
			following inequalities:
			\begin{align}\label{Nasum1}
				\frac{2\alpha Da}{c_s}>\left[\frac{k_L\sqrt{c_p}\alpha_4}{\alpha_3}+{\sqrt{2}k_0}
				(\sqrt{|\Omega|}+\sqrt{c_p}\alpha_7)\right],~~2\alpha_2\geq\frac{\phi^2_s}{a_0},\\
				\label{Nasum2}
				\frac{4\alpha_1Da}{c_kc_s}>\left[\frac{3k_L\sqrt{c_p}\alpha_4}{\alpha_3}
				+{\sqrt{2}k_0}(\sqrt{|\Omega|}+\sqrt{c_p}\alpha_7)\right],
			\end{align}
			where \(\alpha_7= \frac{1}{{\left(\frac{2\alpha_1}{c_k}-\frac{\sqrt{2c_p}\,c_sk_0\alpha_4}{\alpha_3Da}\right)}}\left[
			\sqrt{c_p}||\mathbf{b}_s||_{0,\Omega}+\frac{\alpha_4}{\alpha_3}\left(\phi_s+\frac{\sqrt{2|\Omega|}\,c_sk_0}{Da}\right)\right].\)
\item      In case (b) that is when \(\boldsymbol{\varsigma}=\nabla\cdot\mathbf{U}_s\), the sub-linear growth condition \eqref{P5} becomes 
      \begin{equation}\label{P59a}
          ||\mathbf{K}(\nabla\cdot\mathbf{U}_s)||_{0,\Omega}\leq\sqrt{2}k_0(\sqrt{|\Omega|}+||\nabla\cdot\mathbf{U}_s||_{0,\Omega}).
      \end{equation}
The specific structure $\mathbf{K}(\nabla\cdot\mathbf{U}_s) = [\gamma_1 + \gamma_2 \nabla\cdot\mathbf{U}_s]\mathbf{I}$ falls under this case, and we can clearly see that $\mathbf{K}$ satisfies \eqref{P59a}. In this case, existence and uniqueness analysis can be developed based on similar lines as in Theorem \ref{THM3}. Indeed, we have the following theorem.
   \end{itemize}
		\end{rem}
		\begin{thm}\label{THM39} 
			Assume that the data and parameters satisfy the assumption
			{\bf (A)}, $(i)$ of (\ref{P1}) and (\ref{P3}), \eqref{INEQ29}, \eqref{INEQ39} and \(\mathbf{K}\) satisfies the growth condition \eqref{P59a}. Further, if the following parameter constraint
			\begin{align}\label{Rassup19}
				\frac{2\alpha_1}{c_k}>\frac{\sqrt{2}\,c_sk_0\alpha_4}{\alpha_3Da}		
			\end{align}
			is satisfied, then  \eqref{WFE1} - \eqref{WFE3} together with \(\boldsymbol{\varsigma}=\nabla\cdot\mathbf{U}_s\)
			has a
			solution $(\mathbf{V}_f, \mathbf{U}_s,P)\in H^1(\Omega)^d \times H^1_0(\Omega)^d \times L^2(\Omega)$
			such that
			\begin{align}
				||(\mathbf{V}_f,P)||_{\mathbb{Y}} \leq\frac{\alpha_4}{\alpha_3},~~~~
				||\nabla\mathbf{U}_s||_{0,\Omega}\leq \alpha_8
			\end{align}
			where $|\Omega|$ denotes the area or volume of the domain $\Omega$ and \[\alpha_8=\frac{1}{{\left(\frac{2\alpha_1}{c_k}-\frac{\sqrt{2}\,c_sk_0\alpha_4}{\alpha_3Da}\right)}}\left[
			\sqrt{c_p}||\mathbf{b}_s||_{0,\Omega}+\frac{\alpha_4}{\alpha_3}\left(\phi_s+\frac{\sqrt{2|\Omega|}\,c_sk_0}{Da}\right)\right].\] 
Further, the solution is unique subject to the following constraints
   \begin{align}\label{Nasum19}
				\frac{2\alpha Da}{c_s}>\left[\frac{k_L\alpha_4}{\alpha_3}+{\sqrt{2}k_0}
				(\sqrt{|\Omega|}+\alpha_8)\right],~~2\alpha_2\geq\frac{\phi^2_s}{a_0},\\
				\label{Nasum29}
				\frac{4\alpha_1Da}{c_kc_s}>\left[\frac{3k_L\alpha_4}{\alpha_3}
				+{\sqrt{2}k_0}(\sqrt{|\Omega|}+\alpha_8)\right].
			\end{align}
		\end{thm}

		{\bf Proof:} The proof of this theorem is similar to Theorem \ref{THM3}, we omit the details.
\begin{table}
		\caption{\small{Dimensionless poro-elasto-hydrodynamics parameters corresponding to tumor tissue with their value range}.} \label{Tab2}
		\begin{minipage}{\textwidth}
			\tabcolsep=8pt
			\begin{tabular}{|c|c|c|}
				\hline\hline
				{Dimensionless parameter}  & {Range of values}
				& {Supporting References}\\ 
				\hline
				$Da$   & $10^{-4}-10^{-1}$ & \cite{dey2016hydrodynamics,dey2018vivo}\\\hline
				$\alpha_{t}$   & $0<\alpha_{t}\leq 10$ & \cite{netti1997macro}\\\hline
				$\varrho_{t}$ & $10^2\leq \varrho_{t}\leq 10^5$ & \cite{dey2018vivo}\\\hline
				$\nu_{p}$ &  $0.45\leq \nu_{p} \leq 0.49$ & \cite{roose2007mathematical,islam2020non}\\\hline
				$\phi_f$   & $0.6\leq \phi_f \leq0.8$ &
				\cite{dey2016hydrodynamics,dey2018vivo}\\\hline 
				$K$  
				& $0.00006\leq K \leq1.4$ &
				\cite{dey2018vivo} \\ 
				\hline\hline
			\end{tabular}
		\end{minipage}
	\end{table} 
		
		\subsection{Comments on parameter restrictions}
		It may be noted that the existence and uniqueness results  (e.g. see Theorems \ref{T1}, \ref{T1ab}, \ref{T10}, \ref{T1ab9}, \ref{THM4}, \ref{THM6}, \ref{THM3} and \ref{THM39}, etc.) that
		are established in this work hold subject to certain parameter restrictions see e.g. \eqref{INEQ2}, \eqref{INEQ3}, \eqref{INEQ29}, \eqref{INEQ39},  \eqref{INEQ4}, \eqref{P4}, \eqref{INEQ5}, \eqref{P14}, \eqref{Rassup1}, \eqref{Nasum1}, \eqref{Nasum2}, \eqref{Rassup19}, \eqref{Nasum19} and \eqref{Nasum29} etc. Such a situation is typical in the case of multiphase mixture models where extra care needs to be paid to the physical admissibility of the parameters. See also Vromans et al. \cite{vromans2019parameter}.
		Further,				some of the parameters in these inequalities do show a dependency
		on material properties of the tissue, like Poisson ratio, Young's
		modulus etc. Thus, one has to ensure that these inequalities
		satisfy simultaneously. This certainly depends on the choice of
		relevant parameters. We have collected data from the existing
		literature on parameters like $Da$, $\alpha_t,$ $\varrho_t$ etc.
		which are relevant to biological tissues.
		We have then verified
		that there do exist parameter combinations within the given ranges in Table \ref{Tab2} that obey all the
		inequalities. This ensures that these assumptions
		\eqref{INEQ2}, \eqref{INEQ3}, \eqref{INEQ29}, \eqref{INEQ39},  \eqref{INEQ4}, \eqref{P4}, \eqref{INEQ5}, \eqref{P14}, \eqref{Rassup1}, \eqref{Nasum1}, \eqref{Nasum2}, \eqref{Rassup19}, \eqref{Nasum19} and \eqref{Nasum29} etc. can be interpreted from the point
		of view of various applications. We list the dimensionless
		poro-elasto-hydrodynamics parameters in \ref{Tab2}.
	Further, we choose $\Omega$ as a d-dimensional ($d=3$) sphere of
	unit radius (in dimension 1mm) then $|\Omega|=\frac{4\pi}{3},$
	$|\partial\Omega|=4\pi.$ We set $\mathbf{b}_i=0,$ $i\in\{f,s\},$
	$\mathbf{T}_\infty=(1,0,0),$ and  choose the numerical values for
	constants
	$c_k(\Omega)>0,~c_p(\Omega)>0,~c_s(d)>0,~c_t(\Omega,d)>0$ as
	follows: $c_k=2 ~(or~3)$ for $H^1_0(\Omega)^d$ (or
	$H^1(\Omega)^d$) \mbox{\cite{bernstein1960korn,salsa2016partial}}
	and $c_p=1/2,$ $c_s=1/2$ \mbox{\cite{attouch2014variational}},
	$c_t=2$ \mbox{\cite{salsa2016partial}}. To verify the
	inequalities \eqref{INEQ2}, \eqref{INEQ3}, \eqref{INEQ29}, \eqref{INEQ39},  \eqref{INEQ4}, \eqref{P4}, \eqref{INEQ5}, \eqref{P14}, \eqref{Rassup1} and \eqref{Rassup19}, we use the following combination of parameters
	from Table \mbox{\ref{Tab2}}: $Da=2\times 10^{-2},$
	$\alpha_t=1,$ $\varrho_t=10^4,$ $\nu_p=0.45,$ $k_1=0.5,$
	$k_2=1.4,$ $\phi_s=0.4,$ $L_rA_r=1,$ $k_L=2\times 10^{-3},$
	 and $0<k_0\leq 1.$
	Further, the
	restrictions \mbox{(\ref{Nasum1})}, \mbox{(\ref{Nasum2})} and \eqref{Nasum19},  \eqref{Nasum29} which
	ensure the uniqueness for arbitrary $\mathbf{K}$ hold when we
	choose $Da=10^{-1},$ and $k_0=10^{-2}$ with the remaining
	parameters as chosen above.
	Note that the choice of parameters mentioned above may not be unique; there could be other parameter combinations for which these restrictions hold true.
 
	\section{Summary}
	
 In this work, we have modeled the poro-elasto-hydrodynamics that mimic an in-vitro solid tumor using biphasic mixture theory. We simplified the generic biphasic mixture equations using certain assumptions based on the biological context and treated hydraulic resistivity as anisotropic, which depends on the deformation. This made our model nonlinear and coupled.
We derive an equivalent variational (or weak) formulation and developed existence and uniqueness results using the Galerkin method, monotone operator theory, and fixed-point theory. The detailed analysis is done by considering two cases:
   (a) $\mathbf{K}(\boldsymbol{\varsigma})=\mathbf{K}(\mathbf{U}_s)$
  (b)  $\mathbf{K}(\nabla\cdot\boldsymbol{\varsigma})=\mathbf{K}(\nabla\cdot\mathbf{U}_s)$
In both cases, we first developed existence and uniqueness analysis for auxiliary linear and semilinear sub-problems using the Galerkin method and monotone operator theory. Then, we convert the corresponding coupled nonlinear problem to the fixed-point problem in both cases. Further, to develop the existence of solutions for the corresponding fixed-point problems, we used the Schauder fixed-point theorem. Uniqueness is assured via the Banach contraction theorem.
For the case where $\mathbf{K}$ is not bounded but satisfies the sub-linear growth condition, we have developed the existence and uniqueness results. Moreover, we have collected certain realistic ranges of parameters involved in the model and ensured that the theoretical restrictions derived by us are compatible with these parameter ranges.

	\section*{Acknowledgement}
We are grateful for the valuable constructive suggestions provided by the anonymous referees, which greatly enhanced the quality of our paper. In 
particular, we greatly appreciate the suggestions on alternative methodology to prove one of the main results. M. Alam, one of the authors, is thankful to Dr. Michael Eden and Dr. Debajyoti Choudhuri for their input regarding the analysis.

	\appendix
	\section{}\label{FSpac}
	{\bf Function spaces and useful results\footnote{see
			\cite{salsa2016partial} for function spaces and preliminaries.}:}
	Let $\Omega$ be a bounded, open subset of $\mathbb{R}^d,$ \(\{d=2,3\}.\)
	$L^2(\Omega)$ is the space of all measurable functions $u$ defined
	on $\Omega$ for which
	\begin{equation}\label{Eq0.1}
		||{u}||_{0,\Omega}=\left(\int_{\Omega}|{u}|^2\,d\Omega\right)^{1/2}<+\infty,
	\end{equation}
	In (\ref{Eq0.1}) $||\cdot||_{0,\Omega}$ defines a norm on
	$L^2(\Omega).$ For any $\mathbf{u}=(u_1,u_2,\ldots,u_d)\in
	L^2(\Omega)^d,$ $||\mathbf{u}||_{0,\Omega}$ is defined as
	\begin{equation}\label{Eq0.2}
		||\mathbf{u}||_{0,\Omega}=\left(\int_{\Omega}\sum_{i=1}^{d}|{u}_i|^2\,d\Omega\right)^{1/2},
	\end{equation}
	and for any element $\mathbf{K}=(K_{ij})_{1\leq i,j\leq d}\in
	(L^2(\Omega))^{d\times d},$ we define the norm of $\mathbf{K}$ as
	\begin{equation}\label{Eq0.3}
		||\mathbf{K}||_{0,\Omega}=\left(\int_{\Omega}\sum_{i=1}^{d}\sum_{j=1}^{d}|K_{ij}|^2\,d\Omega\right)^{1/2}.
	\end{equation}
	The symbols $(~,~)_{\Omega},$ and $(~,~)_{\partial\Omega}$ denote  inner products in $L^2(\Omega),$ $L^2(\Omega)^d,$ and $(L^2(\Omega))^{d\times d}$ and in the corresponding trace spaces $L^2(\partial\Omega),$ $L^2(\partial\Omega)^d,$ and $(L^2(\partial\Omega))^{d\times d},$ respectively.\\
	\noindent For any two functions $\mathbf{u}$ and $\mathbf{v}$ the
	inner products $(~,~)_{\Omega},$ and $(~,~)_{\partial\Omega}$ are
	defined as
	$$(\mathbf{u},\mathbf{v})_{\Omega}=\int_{\Omega}\mathbf{u}\cdot\mathbf{v}\,d\Omega,~~ (\mathbf{u},\mathbf{v})_{\partial\Omega}=\int_{\partial\Omega}\mathbf{u}\cdot\mathbf{v}\,d\boldsymbol{\sigma}.$$
	The first-order Sobolev space is defined as\\
	$H^1(\Omega)^d=\{\mathbf{u}\in L^2(\Omega)^d|\nabla{\mathbf{u}}\in
	(L^2(\Omega))^{d\times d}\}$ and the norm of a function
	$\mathbf{u}\in H^1(\Omega)^d$ is defined as
	\begin{equation}\label{Eq0.4}
		||\mathbf{u}||_{1,\Omega}=\left(||\mathbf{u}||^2_{0,\Omega}+
		||\nabla\mathbf{u}||^2_{0,\Omega}\right)^{1/2}.
	\end{equation}
	$H^1_0(\Omega)^d$ denotes the space of functions in $H^1(\Omega)^d$ with zero trace. 
	The dual space of $H^1_0(\Omega)^d$ is denoted by $H^{-1}(\Omega)^d.$
	
	\noindent
	To show the existence of a solution,  we rely on the following results
	\begin{lem}\label{l2}(p.164 \cite{temam2001navier})
		Let $\mathbb{X}$ be a finite-dimensional Hilbert space with scalar
		product $\langle\cdot,\cdot\rangle$ and norm $||\cdot||,$ and let
		${G}$ be a continuous mapping from $\mathbb{X}$ into itself such
		that $$\langle G(x),x\rangle>0~~ \mbox{for}~~ ||x||=r_0>0.$$ Then
		there exists $x\in\mathbb{X},$ with $||x||\leq r_0$ such that
		$$\langle G(x),x\rangle=0.$$
	\end{lem}
	\begin{thm}\label{Schauder} (Schauder's (see p.417 \cite{salsa2016partial}))
		Let \(X\) be a Banach space. Assume that:
		\begin{description}
			\item[(i)] \(A\subset X\) is closed and convex.
			\item[(ii)] \(T: A\rightarrow A\) is continuous.
			\item[(iii)] \(\overline{T(A)}\) is compact in \(X.\)
		\end{description}
		Then \(T\) has a fixed point \(x^*\in A.\)   
	\end{thm}
	\begin{thm}\label{Minty} (Browder-Minty (see p.557 \cite{zeidlernonlinear2}))
		Let \(X\) be a real, separable, reflexive Banach space and let \(T: X\rightarrow X^*  \) (the dual of \(X\)) be bounded, continuous, strongly monotone. Then 
		\[T(u)=g\]
		has a unique solution for each \(g\in X^*.\)
		
	\end{thm}

	\section*{Conflict of interest}	
	The authors declare that there is no conflict of interest.

	\label{lastpage}
	

\begin{thebibliography}{9}
		
		
		\bibitem{dey2016hydrodynamics}
		\textsc{Dey, B. \&  Raja Sekhar, G.~P.}, (2016)\textit{ Hydrodynamics and
			convection enhanced macromolecular fluid transport in soft biological tissues: Application to solid tumor}, Journal of Theoretical Biology, {\bf 395} 62--86.
		\bibitem{byrne1999using}
		{\sc Byrne, H.} (1999) \textit{ Using mathematics to study solid tumour
			growth}. In
		Proceedings of the 9th General Meetings of European Women in Mathematics, New
		York: Hindawi Publishing, 81--107.
		\bibitem{shelton2011mechanistic}
		{\sc  Shelton, S.~E.} (2011) \textit{Mechanistic modeling of cancer tumor
			growth using a  porous media approach}.  Master Thesis, University of North Carolina.
		\bibitem{araujo2004history}
		{\sc Araujo, R.~P. \&  McElwain, D.~S.} (2004) \textit{ A history of the study	of solid tumour growth: the contribution of mathematical modelling}. Bulletin of
		Mathematical Biology, {\bf 66}, 1039--1091.
		
		\bibitem{preziosi2009multiphase}
		{\sc Preziosi, L. \& Tosin, A.} (2009) \textit{Multiphase modelling of
			tumour growth and
			extracellular matrix interaction: mathematical tools and applications}.
		Journal of Mathematical Biology, \textbf{58}, 625--656.
		
		\bibitem{please1998new}
		{\sc Please, C. P., Pettet, G., \& McElwain, D. L. S.} (1998) \textit{ A new approach	to modelling
			the formation of necrotic regions in tumours}, Applied Mathematics Letters,
		\textbf{11}, 89--94.
		
		\bibitem{ambrosi2002closure}
		{\sc Ambrosi, D. and Preziosi, L.} (2002) \textit{ On the closure of mass
			balance models
			for tumor growth}. Mathematical Models and Methods in Applied Sciences, \textbf{12}, 737--754.
		\bibitem{byrne2003modelling}
		{\sc Byrne, H. \& Preziosi, L.} (2003)  \textit{Modelling solid tumour growth
			using the
			theory of mixtures} Mathematical Medicine and Biology, \textbf{20},
		341--366.
		\bibitem{sumets2015boundary}
		{\sc Sumets, P. P., Cater, J. E., Long, D. S., \& Clarke, R. J. } (2015) \textit{ A boundary-integral	representation for biphasic mixture theory, with application to the post-capillary glycocalyx}. Proceedings of the Royal Society A: Mathematical, Physical and Engineering Sciences, {\bf 471}, 20140955.
		\bibitem{barbeiro2010priori}
		{\sc Barbeiro, S. \&  Wheeler, M.~F.} (2010) \textit{ A priori error estimates
			for the
			numerical solution of a coupled geomechanics and reservoir flow model with
			stress-dependent permeability}. Computational Geosciences, \textbf{14},
		755--768.
		\bibitem{jager2011homogenization}
		\textsc{J{\"a}ger, W. and Mikeli{\'c}, A. \& Neuss-Radu, M.} (2011) \textit{Homogenization limit of a model system for interaction of flow, chemical reactions, and mechanics in cell tissues}.
		{SIAM Journal on Mathematical Analysis}, \textbf{43(3)}, 1390--1435.
		\bibitem{cao2014steady}
		{\sc Cao, Y., Chen, S., \& Meir, A. J.} (2014) \textit{ Steady flow in a	deformable porous
			medium}. Mathematical Methods in the Applied Sciences, \textbf{37},
		1029--1041.
		
		\bibitem{alam2018mathematical}
		\textsc{ Alam, M., Dey, B. \&  Raja Sekhar, G.~P.} (2018) \textit{ Mathematical	analysis of
			hydrodynamics and tissue deformation inside an isolated solid tumor}.
		Theoretical and Applied Mechanics, \textbf{45}, 253--278.
		\bibitem{alam2019mathematical}
		\textsc{ Alam, M., Dey, B. \&  Raja Sekhar, G.~P.} (2019) \textit{ Mathematical	modeling and	analysis of hydroelastodynamics inside a solid tumor containing deformable
			tissue}. ZAMM-Journal of Applied Mathematics and Mechanics/Zeitschrift
		f{\"u}r Angewandte Mathematik und Mechanik, e201800223.
		
		\bibitem{alam2021existence}
		\textsc{ Alam, M., Byrne, H. \&  Raja Sekhar, G.~P.} (2021) \textit{ Existence and uniqueness results on biphasic mixture model for an in-vivo tumor}.
		Applicable Analysis, {\bf 101 (15)}, 5442--5468.
		
		\bibitem{barry1990comparison} 
		{\sc Barry, S. \& Aldis, G.}, (1990) \textit{ Comparison of models for flow
			induced
			deformation of soft biological tissue}. Journal of Biomechanics, \textbf{23},
		647--654.
		\bibitem{holmes1990nonlinear}
		{\sc Holmes, M. \& Mow, V.} (1990) \textit{ The nonlinear characteristics of
			soft gels and	hydrated connective tissues in ultrafiltration}. Journal of Biomechanics, \textbf{23}, 1145--1156.
		\bibitem{ateshian2010anisotropic}
		{\sc  Ateshian\, G.~A. \&  Weiss, J.~A.} (2010) \textit{ Anisotropic hydraulic
			permeability
			under finite deformation}. Journal of Biomechanical Engineering, \textbf{132},
		111004.
		
		\bibitem{federico2008anisotropy}
		{\sc Federico, S. and Herzog, W.} (2008) \textit{ On the anisotropy and
			inhomogeneity of
			permeability in articular cartilage}. Biomechanics and Modeling in
		Mechanobiology, \textbf{7}, 367--378.
		
		\bibitem{giverso2019influence}
		\textsc{Giverso, C., \& Preziosi, L.} (2019)
		\textit{ Influence of the mechanical properties of the necrotic core on the
			growth and remodelling of tumour spheroids}.
		{International Journal of Non-Linear Mechanics}, \textbf{108}
		20--32.
		
		\bibitem{barry1991fluid}
		{\sc Barry, S. I., Parkerf, K. H., \& Aldis, G. K.} (1991) \textit{ Fluid flow over a	thin	deformable porous layer}. Zeitschrift f{\"u}r angewandte Mathematik und
		Physik ZAMP, \textbf{42}, 633--648.	
		
		\bibitem{netti1997macro}
		{\sc Netti, P. A., Baxter, L. T., Boucher, Y., Skalak, R., \& Jain, R. K.} (1997) \textit{
			Macro-and microscopic fluid transport in living tissues: Application to solid
			tumors}. AIChE Journal, \textbf{43}, 818--834.
		
		\bibitem{wang1995effect}
		\textsc{ Wang, W. \& Parker, K.~H.} (1995) \textit{ The effect of deformable
			porous surface
			layers on the motion of a sphere in a narrow cylindrical tube}. Journal of
		Fluid Mechanics, {\bf 283}, 287--305.
		\bibitem{sun2002combined}
		\textsc{Sun, S., Riviere, B., \& Wheeler, M. F.} (2002)
		\textit{A combined mixed finite element and discontinuous Galerkin method for miscible displacement problem in porous media}.
		{Recent Progress in Computational and Applied PDEs}, Springer,
		323--351.
		\bibitem{damiano1996axisymmetric}
		{\sc Damiano, E. R., Duling, B. R., Ley, K., \& Skalak, T. C.} (1996) \textit{
			Axisymmetric	pressure-driven flow of rigid pellets through a cylindrical tube lined with a
			deformable porous wall layer}. Journal of Fluid Mechanics, {\bf 314},
		163--189.
		\bibitem{gad1995technical}
		{\sc Gad-el Hak, M.} (1995) \textit{ Stokes hypothesis for a Newtonian,
			isotropic fluid},
		Journal of Fluids Engineering, \textbf{117(1)}, 3-5.
		\bibitem{rajagopal2013new}
		{\sc Rajagopal, K.} (2013) \textit{ A new development and interpretation of
			the
			Navier-Stokes fluid which reveals why the Stokes assumption is inapt},
		International Journal of Non-Linear Mechanics, {\bf 50}, 141--151.
		\bibitem{ccecsmeliouglu2012existence}
		{\sc {\c{C}}e{\c{s}}melio{\u{g}}lu, A. and Rivi{\`e}re, B.} (2012) 
		\textit{Existence of a
			weak solution for the fully coupled Navier--Stokes/Darcy-transport problem}.
		Journal of Differential Equations, {\bf 252}, 4138--4175.
		
		\bibitem{vromans2019parameter}
		\textsc{ Vromans, A. J., van de Ven, A. A. F., \& Muntean, A.} (2019) \textit{ Parameter delimitation of the weak solvability for a pseudo-parabolic system coupling chemical reactions, diffusion and momentum equations}. Adv Math Sci Appl, \textbf{28}, 273--311.
		
		\bibitem{dey2018vivo}
		{\sc Dey, B., Raja Sekhar, G. P., \& Mukhopadhyay, S. K.} (2018) \textit{ In vivo mimicking model for solid tumor
			towards hydromechanics of tissue deformation and creation of necrosis}.
		{Journal of Biological Physics}, {\bf 44(3)}
		361--400.
		\bibitem{roose2007mathematical}
		{\sc Roose, T., Chapman, S. J., \& Maini, P. K.} (2007) \textit{ Mathematical
			models of avascular tumor growth}. SIAM Review, {\bf 49}, 179--208.
		
		\bibitem{islam2020non}
		\textsc{Islam, M., Tang, S., Liverani, C., Saha, S., Tasciotti, E., \& Righetti, R.} (2020)
		\textit{Non-invasive imaging of Young's modulus and Poisson's ratio in cancers in vivo}. {Science Reports}, \textbf{10}, 1-12.
		\bibitem{bernstein1960korn}
		{\sc Bernstein, B. \& Toupin, R.} (1960) \textit{Korn inequalities for the
			sphere and circle}. Archive for Rational Mechanics and Analysis,
		{\bf 6}(1), 51--64.
		\bibitem{salsa2016partial}
		{\sc Salsa, S.} (2016) \textit{Partial Differential Equations in Action:
			From Modelling	to Theory}. {\bf 99}, Springer.
		\bibitem{attouch2014variational}
		{\sc Attouch, H., Buttazzo, G., \& Michaille, G.} (2014) \textit{Variational analysis in Sobolev and BV spaces:
			applications to PDEs and optimization}. SIAM.
		\bibitem{temam2001navier}
		{\sc Temam, R.} (2001) \textit{ Navier-Stokes Equations: Theory and Numerical
			Analysis}.
		{\bf 343}, American Mathematical Society.
		
		\bibitem{zeidlernonlinear2}
		{\sc  Zeidler, E.} (2013) \textit{Nonlinear Functional Analysis and Its Applications: II/B: Nonlinear Monotone
			Operators}.
		Springer Science \& Business Media.
		
		\bibitem{camano2016augmented}
		{\sc  Camano, J.,  Gatica, G. N., Oyarz{\'u}a, R. \& Tierra, G.} (2016) \textit{An Augmented Mixed Finite Element Method for the Navier--Stokes Equations with Variable Viscosity}.
		{SIAM Journal on Numerical Analysis}, {\bf 54(2)}, 1069--1092.
		
		\bibitem{caucao2020fully}
		{\sc  Caucao, S.,  Gatica, G. N., Oyarz{\'u}a, R. \& S{\'a}nchez, N.} (2020) \textit{A fully-mixed formulation for the steady double-diffusive convection system based upon Brinkman--Forchheimer equations}.
		{Journal of Scientific Computing}, {\bf 85(2)}, 37.
		
		
		
		
		
		
		
		
		
	\end{thebibliography}
\end{document}